\documentclass{gtpart}     
\gtart  

\usepackage{geometry}              
\usepackage{setspace}
\usepackage{bm}
\usepackage{color}
\usepackage{tikz}
\usetikzlibrary{arrows,decorations.markings,calc,shapes}
\title[Untwisted automorphisms of RAAGs]{Outer Space for Untwisted Automorphisms of Right-angled Artin Groups}

\author{R. Charney}
\givenname{Ruth}
\surname{Charney}
\address{Brandeis University, Waltham, MA 02453}
\email{charney@brandeis.edu}
\urladdr{http://people.brandeis.edu/~charney/}

\author{N. Stambaugh}
\givenname{Nate}
\surname{Stambaugh}
\address{Western Governors University, Salt Lake City, UT 84107 }
\email{nate.stambaugh@wgu.edu}
\urladdr{}%

\author{K. Vogtmann}
\givenname{Karen}
\surname{Vogtmann}
\address{University of Warwick, Coventry CV4 7AL, United Kingdom}
\email{k.vogtmann@warwick.ac.uk}
\urladdr{https://www2.warwick.ac.uk/fac/sci/maths/people/staff/karen\_vogtmann/}%


\arxivreference{math.GR/1212.4791}
\arxivpassword{9gveb}

\volumenumber{}
\issuenumber{}
\publicationyear{}
\papernumber{}
\startpage{}
\endpage{}
\doi{}
\MR{}
\Zbl{}
\received{}
\revised{}
\accepted{}
\published{}
\publishedonline{}
\proposed{}
\seconded{}
\corresponding{}
\editor{}
\version{}

\makeatletter
\newtheorem*{rep@theorem}{\rep@title}
\newcommand{\newreptheorem}[2]{%
\newenvironment{rep#1}[1]{%
 \def\rep@title{#2 \ref{##1}}%
 \begin{rep@theorem}}%
 {\end{rep@theorem}}}
\makeatother

\newtheorem{theorem}{Theorem}[section]
\newreptheorem{theorem}{Theorem} 
\newtheorem*{corollary*}{Corollary}
\newtheorem{corollary}[theorem]{Corollary}
\newtheorem{lemma}[theorem]{Lemma}
\newtheorem{proposition}[theorem]{Proposition}  
\theoremstyle{definition}
\newtheorem{definition}[theorem]{Definition}
\newtheorem{remark}[theorem]{Remark}
\newtheorem{example}[theorem]{Example}
\newtheorem{examples}[theorem]{Examples}
\newtheorem*{notation}{Notation}
\newcommand{\midarrow}{\draw[postaction={decorate}]}

\newcommand\plainrose[2]{ \node at ({#1},{#2}) [dot]  {}; 
                                         \draw [thick] ({#1},{#2}) to [out=135,in=90] ({#1-1},{#2})  to [out=270,in=225] ({#1},{#2});
                                         \draw  [thick] ({#1},{#2}) to [out=-45,in=270] ({#1+1},{#2}); \draw [thick] ({#1+1},{#2}) to [out=90,in=45] ({#1},{#2})}
 \newcommand\dottedrose[2]{ \node at ({#1},{#2}) [dot]  {}; 
                                         \draw [thick,dotted] ({#1},{#2}) to [out=135,in=90] ({#1-1},{#2});
                                         \draw [thick]  ({#1-1},{#2}) to [out=270,in=225] ({#1},{#2});
                                         \draw  [thick] ({#1},{#2}) to [out=-45,in=270] ({#1+1},{#2});
                                          \draw [thick,,dotted] ({#1+1},{#2}) to [out=90,in=45] ({#1},{#2})}
\newcommand{\iso}{\cong}
\newcommand{\G}{\Gamma}
\newcommand{\GW}{$\G$-Whitehead }  %
\newcommand{\AG}{{\rm A}_\G}            
\newcommand{\Out}{{\rm Out}(\AG)}     
\newcommand{\TAG}{{\rm T}(\AG)}      
 \newcommand{\UAG}{{\rm U}(\AG)}     
 \newcommand{\OG}{\ensuremath{{\mathcal{O}(\AG)}}} 
\newcommand{\KG}{{\rm K}_\Gamma} 
 \newcommand{\Hcal}{{\mathcal H}}      
\newcommand{\Kcal}{{\mathcal K}}       
\newcommand{\Pbar}{{\overline P}}      
\newcommand{\Pbold}{{ \mathbf P}}      
\newcommand{\PP}{{\boldsymbol{\Pi}} } 
\newcommand{\Qbold}{{\mathbf Q}}        
\newcommand{\Rbold}{{\mathbf R}}         
\newcommand{\Mhat}{{\mathbf M}}           
 \newcommand{\SG}{\mathbb S}                 
 \newcommand{\SGP}{{\mathbb S^{\Pbold}}} 
  \newcommand{\SGPP}{{\mathbb S^{\PP}}}  
 \newcommand{\rose}{marked Salvetti} 
 \newcommand{\roses}{marked Salvettis} 
\newcommand{\Xbold}{{\mathbf X}}
\newcommand{\Ybold}{{\mathbf Y}}
\newcommand{\Gcal}{{\mathcal G}} 
\newcommand{\word}{{\bf w}}           
\newcommand{\St}{\mathfrak g}       
\newcommand{\0}{\mathbf 0}
\newcommand{\vcd}{\textsc{vcd}}

\begin{document}

 \begin{abstract} For a right-angled Artin group $\AG,$ the untwisted outer automorphism group $\UAG$ is the subgroup of $\Out$  generated by all of the Laurence-Servatius generators except  twists (where a {\em twist} is an automorphisms of the form $v\mapsto vw$ with $vw=wv$).  We  define a space $\Sigma_\G$ on which $\UAG$ acts properly and prove that $\Sigma_\G$ is contractible, providing a geometric model for $\UAG$ and its subgroups.    We also propose a geometric model for all of $\Out$ defined by allowing more general markings and metrics on points of $\Sigma_\G$.  
  \end{abstract}

\maketitle


\section{Introduction}

A free group is defined by giving a set of generators with no relations; in particular, none of the generators commute.  
A free abelian group is defined by giving a set of generators which all  commute, and no other relations.  
Finitely-generated free and free abelian groups are  examples of {\it right-angled Artin groups}  (RAAGs for short): a general RAAG is  defined by giving a finite set of generators, {\it some} of which commute, and no other relations.  
A convenient way of describing a RAAG is by drawing a graph $\G$ with one vertex for each generator and an edge between each pair of commuting generators; the resulting RAAG is denoted $\AG$.  RAAGs and their subgroups are important sources of examples and counterexamples in geometric group theory (see \cite{Cha} for a survey) and have recently played a key role in the solution of Thurston's conjectures on the structure of hyperbolic 3-manifolds (see ~\cite{AGM}).  

Automorphism groups of RAAGs have received  less attention, with the notable exception of  $\AG=\mathbb Z^n$ and $\AG=F_n$.  Since it is easy to determine the center of any $\AG$ the inner automorphisms of $\AG$ are well-understood, so it remains to study the outer automorphism group $\Out$.   The groups ${\rm Out}(F_n)$ and ${\rm Out}(\Z^n)={\rm GL}(n,\Z)$ have been shown to have many features in common, and it is natural to ask whether these features are in fact shared by all $\Out$. On the other hand, there are important differences between ${\rm GL}(n,\Z)$ and ${\rm Out}(F_n)$ (such as the fact that ${\rm Out}(F_n)$ is not linear!), so we are also interested in how the structure of $\G$ affects the group-theoretic properties of $\Out$. 

In previous work we have explored properties of $\Out$ using inductive local-to-global ideas,  based ultimately on the fact that an outer automorphism of $\AG$ must send certain special subgroups ${\rm A}_{{\rm st}[v]}$ to conjugates of themselves \cite{CCV,ChVo,BCV, ChVo2}.   In this paper we take a more uniformly global approach by introducing an ``Outer space" for any $\AG$ which should play the role of the symmetric space $Q_n={\rm SL}(n,\mathbb R)/{\rm SO}(n)$ in the study of ${\rm SL}(n,\Z)$ or of Outer space $\mathcal O_n$ in the study of ${\rm Out}(F_n)$, i.e. it  should be a contractible finite-dimensional space with a proper action of  $\Out$.   
 
Among the many possible ways of defining  $Q_n$ and $\mathcal O_n$ are   as spaces of free cocompact actions  (of $\Z^n$ on $\mathbb R^n$ or of $F_n$ on simplicial trees) or as spaces of marked metric spaces (flat tori with fundamental group $\Z^n$ or metric graphs with fundamental group $F_n$, equipped with a homotopy equivalence called a {\em marking} to a fixed torus or graph).   All $\AG$ act freely and cocompactly on CAT(0) cube complexes, so it is natural to try to define Outer space in general in the first way, as a space of actions. In \cite{CCV} we  were motivated by this idea but were unable to prove contractiblity of any appropriate space of CAT(0) actions; instead we looked at local data one would obtain from such an action and defined a point of Outer space to be such a data set (whether or not it actually comes from an action).    This trick was successful for  RAAGs defined by graphs which contain no triangles (called 2-dimensional RAAGs), but the methods do not generalize to higher dimension.  

In this paper we take the second approach, defining a space whose points are  marked metric spaces. For every RAAG $\AG$ there is a standard minimal  non-positively curved cube complex with fundamental group $\AG$, called the {\em Salvetti complex}.   We build our space out of somewhat more general cube complexes called {\em $\G$-complexes}; these are homotopy equivalent to Salvetti complexes via an elementary operation called {\em hyperplane collapse}.    For $\AG=\Z^n$ the only $\G$-complex is the Salvetti complex, which is an $n$-torus.  For $\AG=F_n$ the Salvetti complex is a rose with $n$ petals, $\G$-complexes are graphs with no univalent or bivalent vertices or separating edges and hyperplane collapse amounts to collapsing a maximal tree.    
  
We restrict attention to the subgroup $\UAG$ of $\Out$ generated by inversions, graph automorphisms, partial conjugations and transvections of the form $v\mapsto vw$, where $v$ and $w$ are generators of $\AG$ with $vw\neq wv$.    By a theorem of Laurence and Servatius these generate $\Out$ unless there are vertices $v$ and $w$ with ${\rm st}(v)\subseteq {\rm st}(w)$, in which case we must also add {\em twists} sending $v\mapsto vw=wv$.    We refer to $\UAG$ the {\em untwisted} subgroup of $\Out$.  One reason for interest in this subgroup, even when it is not equal to all of $\Out$, is that   for any $\G$  the kernel ${\rm I}\AG$ of the natural map $\Out\to {\rm GL}(n,\Z)$  is contained in $\UAG$.

 We build a simplicial complex $\KG$ on which   $\UAG$   acts properly and cocompactly, so that $\KG$ is quasi-isometric to $\UAG$ by the Schwarz-Milnor lemma. The vertices of $\KG$ are marked $\Gamma$-complexes, where a marking is an ``untwisted" homotopy equivalence to a fixed Salvetti complex.  Two vertices are joined by an edge if one is obtained from the other by a hyperplane collapse, and the resulting graph is then completed to a flag complex to form $\KG$.  Our main theorem is 
 \begin{reptheorem}{contract} For any right-angled Artin group $\AG$, the  complex $\KG$ is contractible.
\end{reptheorem}
As an immediate corollary we obtain
\begin{corollary*}   The dimension of $\KG$ is an upper bound for the virtual cohomological dimension of the untwisted subgroup $\UAG$ of $\Out$.
\end{corollary*}
 The dimension of $\KG$ is always finite.   In Section~\ref{examples} we compute this dimension   for a few examples of $2$-dimensional RAAGs.   Upper and lower bounds for the $\vcd$ of $\Out$ when $\AG$ is $2$-dimensional were computed in \cite{BCV} but these bounds  seldom agree.  Our computations of the dimension of $\KG$ improve on the bounds given in [3] and in some cases,
determine the exact $\vcd$ of $\UAG.$ 

We can obtain a larger space $\Sigma_\G$ by equipping  $\G$-complexes  with metrics in which the cubes are rectilinear parallelepipeds.  
The simplicial complex $\KG$ then embeds as an equivariant deformation retract of $\Sigma_\G$.  
For a free group, $\UAG=\Out$, $\Sigma_\G$ is the usual (reduced) Outer space $\mathcal O_n$, and $\KG$ is the spine of Outer space. For a free abelian group of rank $n$,  $\UAG$ is the finite group of signed permutation matrices,  $\Sigma_\G$ is the positive orthant in $\R^{n-1},$ and $\KG$ is a single point.  

To get a space on which the entire group $\Out$ acts, we allow the markings to be arbitrary homotopy equivalences and the metrics on $\G$-complexes to be ``twisted" in certain restricted ways.  This is discussed briefly in the last section of the paper.  

Our description of  marked $\G$-complexes and proof of  contractibility  of $\KG$ are modeled on Culler and Vogtmann's original  proof that Outer space for a free group is contractible \cite{CuVo}. The idea is that  $\KG$ is the union of the simplicial stars of the \roses, and we assemble all of $\KG$ by attaching these stars one at a time, making sure that at each stage we are gluing along a contractible subcomplex.   The order in which we attach the stars is  determined by a Morse function which measures the lengths of conjugacy classes of $\AG$ under the marking of the Salvetti.  The proof that the subcomplexes along which we glue are contractible requires understanding how this Morse function changes under basic automorphisms; this depends on a generalization of the classical Peak Reduction algorithm for free groups.  A version of Peak Reduction for RAAGs was established by M.~Day in \cite{Day,Day2}.  We require a stronger version (see Theorems~\ref{peak} and ~\ref{s-peak}) and give an independent proof.

We will make use of the standard notions of non-positively curved cube complexes and hyperplanes, and we refer the reader to \cite{HaWi} for these concepts.   

Charney and Vogtmann would like to thank the Mittag-Leffler Institute in Stockholm and the Forschungsinstitut f\"ur Mathematik in Zurich for their hospitality during the development of this paper.  Charney was partially supported by NSF grant DMS-1106726.
    
\section{Whitehead automorphisms}

In this section we recall some basic facts about right-angled Artin groups and their automorphisms.  Fix a right-angled Artin group   $\AG$   with generating set $V$ = vertices$(\G)$.  

\subsection{Partial orders}
Recall from \cite{ChVo} that the relation ${\rm lk}(v)\subseteq {\rm st}(w)$ for $v,w$ vertices of $\Gamma$ is denoted $v\leq w$.  Vertices are called {\em equivalent} if $v\leq w$ and $w\leq v$, and we write $v\sim w$.  (The justification for this notation is that $\leq$ is a partial order on equivalence classes of vertices.)  If $v$ is adjacent to $w$ then $v\leq w$ if and only if ${\rm st}(v)\subseteq {\rm st}(w)$, and if $v$ is not adjacent to $w$ then $v\leq w$ if and only if  ${\rm lk}(v)\subseteq {\rm lk}(w)$. 
 
When considered as elements of $\AG$, each element  $v\in V$ has an inverse $v^{-1}$, and we will often work with the symmetric set $V^\pm = \{v,v^{-1}\mid v\in V\}$.  For $x,y\in V^{\pm}$ we say $x\leq y$ etc.
if the corresponding vertices of $\Gamma$ satisfy the relation.

\subsection{Generators for \texorpdfstring{$\Out$}{Out(Fn)}} Laurence and Servatius (\cite{Lau, Ser}) proved that the following simple types of automorphisms generate all of ${\rm Aut}(\AG)$ (and hence their images generate $\Out$):
\begin{enumerate}
\item  An  automorphism of  the graph $\G$ permutes the vertices $V$ and induces an automorphism of $\AG$, called a {\em graph automorphism}.

\item If $v\in V$, the map sending $v\mapsto v^{-1}$ and fixing all other generators is an automorphism of $\AG$, called an {\it inversion}.

\item If $v\leq w$, then the map sending $v\mapsto vw$ and fixing all other generators is an automorphism of $\AG$, called a {\em transvection}.  If $v$ is adjacent to $w$ this is called a {\em twist}, and if not it is called a {\em fold}.   

\item If $C$ is a component of $\G\setminus {\rm st}(v)$, then the map sending $x\to vxv^{-1}$ for every vertex $x$  of $C$ and fixing all other generators is   an automorphism of $\AG$, called a {\em partial conjugation}.  

\end{enumerate}

Twists play a specialized role in the study of $\Out$, and we define $\TAG$ to be the subgroup of $\Out$ generated by these.  We then define  $\UAG$ to be the subgroup generated by all other types of generators, i.e.  graph automorphisms, inversions, partial conjugations and folds.  In the terminology of \cite{Day}, elements of $\UAG$ are called \emph{long-range} automorphisms.  The subgroup $\UAG$ will be the focus of the rest of the paper.

\subsection{\texorpdfstring{$\boldsymbol\G$}{G}-Whitehead  partitions} There is a larger generating set for $\UAG$ which is more natural for our purposes.  This larger set includes simple combinations of folds and partial conjugations; it consists of automorphisms of $\AG$ which are induced by Whitehead automorphisms of the free group ${\rm F}(V)$.  We recall that a Whitehead automorphism of ${\rm F}(V)$ is determined by a pair $(P,m)$, where $P\subset V^\pm$ has at least 2 elements, and $m\in P$ with $m^{-1}\not\in P$.  The automorphism $\varphi=(P,m)$  is given by
$$ \varphi(v) = \begin{cases} m^{-1}&\textrm{if $v=m$}\cr
       vm^{-1} &\textrm{if  $ v \in P$ and $v^{-1} \notin  P$}\cr
	mv &\textrm{if  $ v^{-1} \in  P$ and $v \not\in  P$}\cr
	mvm^{-1} &\textrm{if  $v, v^{-1} \in  P$}\cr
	v &\textrm{otherwise}
	\end{cases} $$
We remark that it is more usual to define a Whitehead automorphism with $\varphi(m)=m$, but as we will see in Lemma~\ref{blowupcollapse} below, setting $\varphi(m)=m^{-1}$ corresponds more naturally with the geometric version of a Whitehead move. (This was originally observed by Hoare in \cite{Hoare}).  With this definition, $\varphi$ is an involution, so $\varphi=\varphi^{-1}$. 
Replacing $P$ by its complement, $P^*$  and $m$ by $m^{-1}$ changes $\varphi$ by an inner automorphism (conjugation by $m$).  Thus, $(P,m)$ and $(P^*,m^{-1})$ determine the same outer automorphism.

Not every Whitehead automorphism of ${\rm F}(V)$ induces an automorphism of $\AG$, and even if it does the induced automorphism may be trivial (e.g. conjugating $v$ by an adjacent $w$). Both of these problems are solved by the following definition.
\begin{definition}\label{Whitehead} Let $P\subset V^{\pm}$ have at least 2 elements, including some $m\in P$ with $m^{-1}\not\in P$.  Then $(P,m)$ is a {\em  \GW  pair} if  
\begin{enumerate}
\item no element of $P$ is adjacent to $m$,
\item if $x\in P$ and $x^{-1}\not\in P$, then $x\leq m$, and
\item if $v, v^{-1}\in P$ then $w,w^{-1}\in P$ for all $w$ in the same component as $v$ of $\G\setminus {\rm st}(m)$.
\end{enumerate}
\end{definition}

\begin{lemma}\label{induced}  If $(P,m)$ is a  \GW  pair, then  the Whitehead automorphism of ${\rm F}(V)$ defined  by $(P,m)$ induces a non-trivial automorphism of $\AG$.  The induced outer automorphism lies in $\UAG$.
\end{lemma}

\begin{proof} Conditions (1)-(3) in the definition of a \GW pair guarantee  that any relation $[v,w]=1$ is preserved by the Whitehead automorphism $(P,m)$.
\end{proof}

\begin{figure}
\begin{center}
 \begin{tikzpicture}[scale=.4]
\draw [densely dotted] (0,0)  ellipse   (4.5cm and .75cm);
\node (lk) at (5.5,.5) {${\rm lk}(m)$};   
\draw  [densely dotted ] (-4,-4.25)  ellipse  (1.25cm and 2cm);  
 \node (C1) at (-4,-5.6) {$C_1$};
\draw  [densely dotted ] (0,-4.25)  ellipse   (1.5cm and 2cm); 
 \node (C2) at (0,-5.6) {$C_2$};
\draw  [densely dotted ] (4,-4.25)  ellipse   (1.25cm and 2cm);
\node (C3) at (4,-5.6) {$C_3$};        
\fill(-1,0) circle  [radius=2pt];
\fill(-3,0) circle  [radius=2pt] node[left] (x) {$x$};
\fill(1,0) circle  [radius=2pt];
\fill(3,0) circle  [radius=2pt];
\fill(-2,3) circle  [radius=2pt] node[above] (m) {$m$};
\fill(0,3) circle  [radius=2pt] node[above] (u) {$u$};
\fill(2,3) circle  [radius=2pt] node[above] (v) {$v$};
\fill(0,-3) circle  [radius=2pt]; 
\fill(-4,-3) circle  [radius=2pt] node [left] (y) {$y$};
\fill(4,-3) circle  [radius=2pt];
\fill(-4,-5) circle  [radius=2pt];
\fill(-1,-5) circle  [radius=2pt];
      \fill(1,-5) circle  [radius=2pt];
        \fill(4,-5) circle  [radius=2pt];
\draw (-4,-5) -- (-4,-3)-- (-3,0) --(-2,3)--(-1,0) -- (-4,-3);
\draw (0,3) -- (-1,0)--(0,-3)--(1,-5) -- (-1,-5)--(0,-3)--(1,0)--(0,3);
\draw (-1,0) -- (1,0)--(-2,3)--(3,0) -- (2,3)--(1,0)--(4,-3)--(3,0);
\draw (4,-5) -- (4,-3);
\end{tikzpicture}
\bigskip
  
   \begin{tikzpicture}[scale=.8] 
                \fill [lightgray!30, rounded corners] (3,1.4) to  ( 9.4,1.4) to (9.4, -1.4) to (4.6,-1.4) to (4.6,.2) to (1.6,.2) to  (1.6,1.4) to (3, 1.4);
                \node (P) at (5.5,1.75) {$P$};
                   \fill [lightgray!30, rounded corners] (.5,1.4) to  (1.4,1.4) to (1.4, -.2) to (4.3,-.2) to (4.3,-1.4) to (-.4,-1.4) to  (-.4,1.4) to (.5, 1.4);
                \node (Pstar) at (.5,1.75) {$P^*$};

   \foreach \x in {0,...,13}
   \fill(\x,1) circle  [radius=2pt]; 
      \foreach \x in {0,...,13}
    \fill(\x,-1) circle  [radius=2pt];
    \node [below] (m) at (2,1) {$m$};
       \node [above] (mi) at (2,-1) {$m^{-1}$};
       \node [below] (u) at (3,1) {$u$};
       \node [above] (ui) at (3,-1) {$u^{-1}$};
              \node [below] (v) at (4,1) {$v$};
       \node [above] (vi) at (4,-1) {$v^{-1}$};
              \node [below] (y) at (5,1) {$y$};
       \node [above] (yi) at (5.2,-1) {$y^{-1}$};
              \node [below] (x) at (10,1) {$x$};
       \node [above] (xi) at (10.2,-1) {$x^{-1}$};
       \draw [densely dotted, rounded corners] (-.25,-1.25) rectangle   (1.25,1.25);
       \node (C3) at (.5,0) {$C_3$};
        \draw [densely dotted, rounded corners] (4.75,-1.25) rectangle   (6.25,1.25);
            \node (C1) at (5.5,0) {$C_1$};
            \draw [densely dotted, rounded corners] (6.75,-1.25) rectangle   (9.25,1.25);
                \node (C2) at (8,0) {$C_2$};
            \draw [densely dotted, rounded corners] (9.75,-1.25) rectangle   (13.25,1.25);
                \node (lkP) at (11.5,0) {${\rm lk}(P)$};  
 \end{tikzpicture} 
\end{center}
\caption{A graph $\G$ and a \GW partition}
\label{partition}
\end{figure}

If $(P,m)$ is a \GW pair, we define 
\begin{align*}
{\rm double}(P) &= \{x\in P \mid \textrm{$x^{-1}$ is also in $P$}\}\\
{\rm single}(P) &=  \{x\in P \mid \textrm  {$x^{-1}$ is not in P}\}\\
{\rm max}(P) &=   \{x \in {\rm single}(P)  \mid x\sim m\}\\
{\rm lk}(P) &= {\rm lk}(m)^\pm
\end{align*}
  
\begin{remark} Here are some elementary observations about \GW pairs.  Since all elements of ${\rm max}(P)$ are equivalent,  ${\rm lk}(P)$ is independent of the choice of $m \in {\rm max}(P)$.  For any $m^\prime\in {\rm max}(P)$, the pair $(P,m^\prime)$ is also a \GW pair.  The set ${\rm max}(P)$ can be recovered from $P$ without reference to $m$ as the set of maximal elements in ${\rm single}(P)$.     Since every $v\in {\rm single}(P)$ is $\leq m$ and is not adjacent to $m$, no two elements of ${\rm single}(P)$ are adjacent to each other.
  
\end{remark}  

By condition (1)  of Definition~\ref{Whitehead} ${\rm lk}(P)$ is disjoint from $P$.  Let  $P^*$ be the complement of  $P \cup {\rm lk}(P)$ in $V^\pm, $   i.e. we have a partition of $V^\pm$ into three disjoint subsets
$$ V^\pm = P + {\rm lk}(P) + P^* $$
It is easy to verify that $(P^*,v^{-1})$ is also a \GW pair for any $v\in {\rm max}(P)$,  that ${\rm lk}(P^*)={\rm lk}(P)$ and that ${\rm max}(P^*)={\rm max}(P)^{-1}$.

\begin{definition} If $(P,m)$ is a \GW pair,  the triple 
$\Pbold= \{P,{\rm lk}(P),P^*\}$ is called a \emph{\GW partition of $V^\pm$},  and $P$ and $P^*$ are called the {\em sides} of $\Pbold.$   
\end{definition}

We remark that a $\G$-Whitehead partition is completely determined by either of its sides.  

\begin{notation}
We will often use $P^\times$ to denote a choice of side of $\Pbold.$  The following notation distinguishes  vertices of $\G$ (as opposed to elements of $V^{\pm}$):
\begin{align*}
{\rm lk}(\Pbold) &= \{v\in V\mid v,v^{-1}\in {\rm lk}(P)\}= \{v\in V\mid v,v^{-1}\in {\rm lk}(P^*)\}\\
{\rm single}(\Pbold)&= \{v\in V \mid \textrm{$v$ and $v^{-1}$ are in different sides of $\Pbold$}\}\\
{\rm double}(\Pbold) &= \{v\in V \mid \textrm{$v$ and $v^{-1}$ are both in the same side of $\Pbold$}\}\\
{\rm max}(\Pbold) &= \{v\in V \mid \textrm{$v$ or $v^{-1}$ is in ${\rm max}(P)$}\}
\end{align*}
\end{notation}

\section{\texorpdfstring{$\boldsymbol\G$}{G}-complexes}
\subsection{Blowing up a single \texorpdfstring{$\boldsymbol\G$}{G}-Whitehead partition} 

\begin{figure}
\begin{center}
\begin{tikzpicture}
[dot/.style={circle, draw, minimum size=1mm,inner sep=0pt,outer sep=0pt, fill}]
\begin{scope}[scale=1.35]
\begin{scope}[decoration={markings,mark = at position 0.5 with {\arrow{stealth}}}] 
\dottedrose{0}{0}; 
\plainrose{0}{2};
\midarrow [thick] (0,0)--(0,2); 
\draw [thin] (-1,0) -- (-1,2);
\draw [thin] (1,0) -- (1,2);
\node (em) [right] at (0,1) {$e_m$};
\node (lkm) [right] at (0,.5) {$S_{{\rm lk}(m)}\times I$};
\begin{scope}[xshift=3cm]
\dottedrose{0}{0}; 
\plainrose{0}{2};
\midarrow [thick] (0,0)--(0,2); 
\draw [thin] (-1,0) -- (-1,2);
\draw [thin] (1,0) -- (1,2);
\node (em) [right] at (0,1) {$e_P$};
\node (lkm) [right] at (0,.5) {$S_{L}\times I$};
\end{scope}
\begin{scope}[xshift=5.5cm]
 \node at (0,0) [dot]  {}; 
 \draw  [thick] (0,0) to [out=-45,in=270] (1,0); \draw [thick, dotted] (1,0) to [out=90,in=45] (0,0);
 \draw  [thick] (0,2) to [out=-45,in=270] (1,2) to [out=90,in=45] (0,2);
\midarrow [thick] (0,0)--(0,2); 
\draw [thin] (1,0) -- (1,2);
\node (em) [right] at (0,1) {$e_n$};
\node (lkm) [right] at (0,.5) {$S_{{\rm lk}(n)}\times I$};
\end{scope}
\begin{scope}[xshift=3cm]
\plainrose{0}{-1};
\draw (-1,-1) .. controls (-1,-2.3) and (1,-2.3) .. (1,-1);
\node (LplusDstar)  at (0,-1.5) {$S_{L+D^*}$};
\end{scope}
\begin{scope}[xshift=3cm]
\dottedrose{0}{3};
 \draw (-1,3) .. controls (-1,4.5) and (1,4.5) .. (1,3);
\draw [thick] (0,3) to [out=45, in=0] (0,4) to [out=180,in=135] (0,3);
 \node (LplusD) [right] at (.2,3.3) {$S_{L+D}$};
\end{scope}
\draw[->,  thin] (3,2.8)--(3,2.2);
\draw[->,  thin] (3,-.8)--(3,-.2);
\draw[->,  thin] (.7,2.4) to [out=60, in = 120] (2.2,2.4);
\draw[->,  thin] (.7,-.4) to [out=-60, in = 240] (2.2,-.4);
\draw[->,  thin] (5.6,-.4) to [out=240, in = -60] (4,-.4);
\draw[->,  thin] (5.6, 2.4) to [out=120, in = 60] (4,2.4);
\end{scope}
\end{scope}
\end{tikzpicture}
\end{center}
\caption{Blowup $\SGP$ of $\SG$}
\label{blowup}
\end{figure}

We begin by recalling the construction of the Salvetti complex $\SG=\mathbb S_\Gamma$.  
Let $n$ be the cardinality of $V$ and let $\mathbb T^n$ denote an $n$-torus with  edges labelled 
$\{e_v \mid v \in V\}$.    Then $\SG$ is the subcomplex of $\mathbb T^n$ consisting of faces whose edges are labelled by mutually commuting sets of vertices.  It is easily verified that $\SG$ is locally CAT(0) (hence aspherical) and has fundamental group $\AG$.  
For any subset $U \subset V$, let $\SG_U$ denote the subcomplex of $\SG$ spanned by the edges labelled $e_u,  u \in U$.  Note that this is isomorphic to the Salvetti complex for the RAAG generated by $U$.

Let $\Pbold=\{P,{\rm lk}(P),P^*\}$ be a \GW partition of $V^\pm$. We  obtain a new cube complex $\SGP$ by {\em blowing up\,}   $\SG$  as follows.
Let $D$ denote the vertices represented in ${\rm double}(P)$, $D^*$ the vertices in ${\rm double}(P^*)$ and $L={\rm lk}(\Pbold)$.  
To construct  $\SGP$:
\begin{itemize}
\item Start with a copy of $\SG_{L}\times [0,1]$.  Label the (unique) 1-cell $e_\Pbold$.
\item Attach a copy of $\SG_{L+ D}$   by identifying the natural subcomplex $\SG_L\subset \SG_{L+ D}$ with $\SG_{L}\times \{1\}$.
\item Attach a copy of $\SG_{L+ D^*}$   by identifying the natural subcomplex $\SG_L\subset \SG_{L+ D^*}$ with $\SG_{L}\times \{0\}$.
\item For each $v\in {\rm single}(\Pbold)$, attach a copy of  $\SG_{{\rm lk}(v)}\times [0,1]$ at its ends using the natural inclusions  $\SG_{{\rm lk}(v)}\times \{1\}\subset \SG_{L}\times \{1\}$ and   $\SG_{{\rm lk}(v)}\times \{0\}\subset \SG_{L}\times \{0\}.$ Label the edge of  $\SG_{{\rm lk}(v)}\times [0,1]$ with $e_v$. Orient it from $0$ to $1$ if $v\in P$ and from $1$ to $0$ if $v^{-1}\in P$.
\end{itemize} 
Figure~\ref{blowup}  may help the reader  visualize this construction.

\begin{remark} If $v$ and $w$ commute, then $\SG$ contains a corresponding torus $T(v,w)$, say with $e_v$ as longitude and $e_w$ as meridian. This torus ``blows up" to the following subcomplex of $\SGP$:

\begin{itemize} 
\item If $v\in {\rm double}(P)$ then $w$ must  be in ${\rm lk}(P)$ or ${\rm double}(P)$.  In either case  $T(v,w)$ gives rise to a torus attached at the vertex of $e_\Pbold$ in $\SG_{L}\times \{1\}$, with $e_v$ as longitude and (the appropriate copy of) $e_w$ as meridian.  
 If $v\in {\rm double}(P^*)$ then the torus is attached at the  vertex of $e_\Pbold$ in $\SG_{L}\times \{0\}$.  
\item If $v\in {\rm single}(\Pbold)$  the longitude of $T(v,w)$ is subdivided into two edges  labeled $e_\Pbold$ and   $e_v$. The meridian loop $w\in {\rm lk}(v)\subset {\rm lk}(P)$ has two representatives, one at each end of $e_\Pbold$, both labeled $e_w$.
 \item If $v$ and $w$ are both in ${\rm lk}(P)$,  then $T$ blows up to the product $T(v,w)\times e_\Pbold\subseteq\SGP$. 

\end{itemize}

\end{remark}

We note the following properties of the blowup:  

\begin{enumerate}
\item $\SGP$ has exactly two vertices which correspond to the two sides of $\Pbold$.
\item  The edges emanating from the $P$-vertex are labelled by the elements of $P \cup {\rm lk}(P)$ plus one extra edge labelled $e_\Pbold$.
Similarly for the $P^*$-vertex.
\item  Two edges at a vertex span a square if and only if they are labelled by commuting generators, or by $e_\Pbold$ and an element of 
${\rm lk}(P)$.
\item The links of the vertices are flag, hence $\SGP$ is non-positively curved.
\end{enumerate}

Note that collapsing the cylinder $\SG_L\times [0,1]$ down to $\SG_L\times \{0\}$ recovers  $\SG$; we call this the {\em canonical collapse $c_\Pbold$}.  For each $m\in {\rm max}(P)$, there is an isomorphism 
$h_m$ of $\SGP$ which interchanges $\SG_L \times e_P$ and $\SG_L \times e_m$.  Let $c_m$ denote the composite map $c_P \circ h_m \colon \SGP \to \SG$.

\begin{lemma}\label{blowupcollapse}  Let $c_\Pbold^{-1}$  be  a homotopy inverse of the canonical collapse.  Then the composition $c_m\circ c_\Pbold^{-1}\colon \SG\to \SGP\to  \SG$ induces the Whitehead automorphism $(P,m)$ on $\AG=\pi_1(\SG)$.       
\end{lemma}

\subsection{Compatible and commuting \texorpdfstring{$\boldsymbol\G$}{G}-Whitehead partitions}

It is possible to build a connected graph with a proper cocompact action of $\UAG$ using just Salvetti complexes and single blowups $\SGP,$ but to make a \emph{contractible} complex we will need to do further blowing up to ``fill in the holes" in this graph. To this end, we make the following definitions.

\begin{definition} Let   $\Pbold=\{P,{\rm lk}(P),P^*\}$ and $\Qbold= \{Q,{\rm lk}(Q),Q^*\}$ be two \GW partitions. 
\begin{enumerate}
\item  Say $\Pbold, \Qbold$ 
\emph{commute} if the equivalence classes of ${\rm max}(P)$ and ${\rm max}(Q)$ are distinct and commute in $\AG$.
\item Say  $\Pbold, \Qbold$ are
\emph{compatible} if either they commute, or $P^\times\cap Q^\times$ is empty for (at least) one choice of sides $P^\times\in \{P,P^*\}  $ and $Q^\times\in \{Q,Q^*\}$.
\end{enumerate}
\end{definition}

\begin{lemma}\label{commute}   Let $\Pbold$ and $\Qbold$ be non-commuting compatible \GW partitions. If $P\cap Q =\emptyset$, then $P\cap {\rm lk}(Q)=\emptyset$, i.e., $P \subset Q^*$ and $Q \subset P^*$.  
\end{lemma}

\begin{proof} Suppose $u\in P\cap {\rm lk}(Q)$, and let $m\in {\rm max}(Q)$.  Then $u \in {\rm lk}(Q)={\rm lk}(m)$  implies that $m\in {\rm lk}(u)$.   If $u\in {\rm single}(P)$  then ${\rm lk}(u)\subseteq {\rm lk}(P)$, so  $m\in {\rm lk}(P)$, contradicting the assumption that $\Pbold$ and $\Qbold$ do not commute.     If $u\in {\rm double}(P)$  then  the fact that $u$ and $m$ are connected by an edge  implies that  either $m \in {\rm lk}(P)$, or $m$ lies in the same component of  $\G -{\rm lk}(P)$ as $u$.      The former contradicts the assumption that $\Pbold$ and $\Qbold$ do not commute,
and the latter implies that $m \in P \cap Q$.   

The last statement follows by symmetry.
 \end{proof}

\begin{remark}\label{maxPQ}   Let $\Pbold$ and $\Qbold$ be non-commuting compatible \GW partitions and suppose that 
$m \in {\rm max}(P) \cap {\rm max}(Q)$.  Then either $P^* \cap Q$ or $P \cap Q^*$ is empty, say $P^* \cap Q = \emptyset$.
Then it follows from the lemma that $Q \subset P$ and setting $R=(P\setminus Q) \cup \{m\},$  a straightforward exercise shows that $(R,m)$  is also a  \GW pair and the corresponding Whitehead automorphisms satisfy
$(P,m) \circ (Q,m)^{-1} = (P,m) \circ (Q,m)= (R,m) \circ i_m$ where $i_m$ is the inversion taking $m \mapsto m^{-1}$.
\end{remark}
 
\begin{lemma}\label{noncommute} Let $\Pbold$ and $\Qbold$ be distinct compatible \GW partitions.  If $\Pbold$ and $\Qbold$ do not commute, then exactly one of $P^\times \cap Q^\times$ is empty.
\end{lemma}

\begin{proof}  Since $\Pbold$ and $\Qbold$ are compatible at least one of the sets $P^\times\cap Q^\times$ is empty, so without loss of generality we may assume $P\cap Q=\emptyset$.  By the previous lemma, it follows that $P\cap Q^*=P$ and $P^* \cap Q=Q$.

Suppose $P^*\cap Q^*=\emptyset$.  Then  any  $m\in {\rm max}(P)$ must have $m^{-1}\in Q$, so ${\rm lk}(P)\subseteq {\rm lk}(Q)$ and similarly ${\rm lk}(Q)\subseteq {\rm lk}(P)$.  Thus $V^{\pm}=P \sqcup{\rm lk}(P) \sqcup  Q$, i.e. $Q=P^*$ and $\Pbold=\Qbold$, contradicting our hypothesis.    
\end{proof}  

 It follows from the lemma that for non-commuting, compatible partitions $\Pbold, \Qbold$ with non-empty intersection   $P^{\times} \cap Q^{\times}$, we can switch sides of 
either $\Pbold$ or $\Qbold$, but not necessarily both, and still get a non-empty intersection.

\subsection{Blowing up compatible collections of \texorpdfstring{$\boldsymbol\G$}{G}-Whitehead partitions}
Now let $$\PP= \{\Pbold_1, \dots ,\Pbold_k\}$$ be a set of pairwise compatible \GW partitions.  We want to simultaneously blow-up $\SG$ along all of the partitions in $\PP$ to obtain a non-positively curved cube complex $\SGPP$.         

The role played by the edge $e_\Pbold$ in the single blowup $\SGP$ will now be played by a cubical subcomplex of a $k$-dimensional cube, $[0,1]^k$. The vertices of this subcomplex will form the vertices of $\SGPP$, and to describe them we make the following definition:
  
\begin{definition}  A \emph{region} of  $\PP$ is a choice of side $P^{\times}_i \in \{P_i, P_i^*\}$ for each $i$ such that
 for $i \neq j$, either $\Pbold_i, \Pbold_j$ commute, or $P^{\times}_i \cap  P^{\times}_j \neq \emptyset$.  
\end{definition}
To each region $R=(P^{\times}_1,\ldots,P^{\times}_k)$ of $\PP$ we associate a vertex $x_R=(a_1,\ldots,a_k)$ of $[0,1]^k$  by the rule 
\begin{equation*}
a_i = 
\begin{cases} 0 & \text{if $P_i^\times = P_i$,}
\\
1 &\text{if $P_i^\times = P_i^*$.}
\end{cases}
\end{equation*}
Let $\mathbb E_\PP$ denote the cubical subcomplex of $[0,1]^k$ spanned by the $x_R$, and label all edges parallel to the $i$-th basis vector with the label $e_{\Pbold_i}$.   

\begin{example} Suppose $\PP=\{\Pbold,\Qbold\}$.  If $\Pbold$ and $\Qbold$ commute, then $\mathbb E_\PP$ is the entire square $[0,1]^2$, with two (parallel) edges labeled $e_\Pbold$ and the other two labeled $e_\Qbold$.    If $\Pbold$ and $\Qbold$ do not commute, then by Lemma~\ref{noncommute} exactly three of $(P,Q), (P^*,Q), (P,Q^*)$ and $(P^*,Q^*)$ are regions, so that $\mathbb E_\PP$ consists of two adjacent edges of the square, one labeled $e_\Pbold$ and one labeled $e_\Qbold$.  

If $\PP$ contains three \GW partitions,  the possibilities for $\mathbb E_\PP$ are illustrated in Figure~\ref{EPQR}.
\end{example}

\begin{figure}
\begin{center}
\begin{tikzpicture}[scale=.5]
 \fill [blue!30]  (0,0) to (3,0) to (4,1) to (4,4) to (1,4) to (0,3) to (0,0);
\draw [blue] (1,1) to (4,1) to (4,4) to (1,4) to (1,1);
 \draw [blue] (0,0) to (1,1) to (4,1) to (3,0) to (0,0);
 \draw [blue] (0,0) to (3,0) to (3,3) to (0,3) to (0,0);
 \draw [blue](0,3) to (1,4) to (4,4) to (3,3) to (0,3);
 \draw [blue](0,0) to (0,3) to (1,4) to (1,1) to (0,0);
 \draw [blue](3,0) to (3,3) to (4,4) to (4,1) to (3,0);
 
\begin{scope}[xshift=5cm]
 \draw [fill=blue!30] (0,0) to (0,3) to (1,4) to (1,1) to (0,0);
  \draw [blue] (0,0) to (0,3) to (1,4) to (1,1) to (0,0);
 \draw [fill=blue!30] (0,0) to (1,1) to (4,1) to (3,0) to (0,0);
  \draw [blue] (0,0) to (1,1) to (4,1) to (3,0) to (0,0);
\draw [very thin, gray] (4,1) to (4,4) to (1,4);
 \draw [very thin, gray] (3,0) to (3,3) to (0,3);
 \draw [very thin, gray](0,3) to (1,4);
 \draw [very thin, gray](3,0) to (3,3) to (4,4) to (4,1) to (3,0);
 \end{scope}
 
 \begin{scope}[xshift=10cm]
 \draw [fill=blue!30](0,0) to (1,1) to (4,1) to (3,0) to (0,0);
  \draw [blue](0,0) to (1,1) to (4,1) to (3,0) to (0,0);
\draw [very thin, gray] (4,1) to (4,4) to (1,4);
 \draw [very thin, gray] (3,0) to (3,3) to (0,3);
\draw [very thin, gray] (0,3) to (1,4) to (4,4) to (3,3) to (0,3);
\draw [very thin, gray] (0,0) to (0,3) to (1,4) to (1,1);
 \draw [very thin, gray]  (3,3) to (4,4) to (4,1);
 \draw [blue!50, line width=.75mm] (0,0) to (0,3);
  \draw(0,0) to (0,3);
\end{scope}

 \begin{scope}[xshift=15cm]
\draw [very thin, gray](1,1) to (4,1) to (4,4) to (1,4) to (1,1);
 \draw [very thin, gray](0,0) to (1,1) to (4,1) to (3,0) to (0,0);
 \draw [very thin, gray](0,0) to (3,0) to (3,3) to (0,3) to (0,0);
 \draw [very thin, gray](0,3) to (1,4) to (4,4) to (3,3) to (0,3);
 \draw [very thin, gray](0,0) to (0,3) to (1,4) to (1,1) to (0,0);
 \draw [very thin, gray](3,0) to (3,3) to (4,4) to (4,1) to (3,0);
  \draw [blue!50, line width=.75mm] (0,0) to (0,3);
    \draw [blue!50, line width=.75mm] (0,0) to (3,0);
      \draw [blue!50, line width=.75mm] (0,0) to (1,1);
      \draw  (0,0) to (0,3);
    \draw   (0,0) to (3,0);
      \draw  (0,0) to (1,1);
\end{scope}

 \begin{scope}[xshift=20cm]
\draw [very thin, gray](1,1) to (4,1) to (4,4) to (1,4) to (1,1);
 \draw [very thin, gray](0,0) to (1,1) to (4,1) to (3,0) to (0,0);
 \draw [very thin, gray](0,0) to (3,0) to (3,3) to (0,3) to (0,0);
 \draw [very thin, gray](0,3) to (1,4) to (4,4) to (3,3) to (0,3);
 \draw[very thin, gray] (0,0) to (0,3) to (1,4) to (1,1) to (0,0);
 \draw[very thin, gray] (3,0) to (3,3) to (4,4) to (4,1) to (3,0);
  \draw [blue!50, line width=.75mm] (3,0) to (3,3);
    \draw [blue!50, line width=.75mm] (0,0) to (3,0);
      \draw [blue!50, line width=.75mm] (0,0) to (1,1);
 \draw (3,0) to (3,3);
    \draw  (0,0) to (3,0);
      \draw  (0,0) to (1,1);
\end{scope}
 \end{tikzpicture}
\end{center}
\caption{Possibilities for $\mathbb E_{\{\Pbold,\Qbold,\mathbf R\}}$}
\label{EPQR}
\end{figure}

The following  lemma guarantees that every set of pairwise compatible \GW partitions has regions associated to it.

\begin{lemma}\label{regions}   Let $\PP= \{\Pbold_1, \dots ,\Pbold_k\}$ be a set of pairwise compatible \GW partitions,  and $(P_1^\times, \dots ,P_k^\times)$  a region for $\PP$.  If a \GW partition $\Pbold$ is compatible with each $\Pbold_i$ then for some choice of sides $P^\times\in  \{P, P^*\}$, $(P_1^\times, \dots ,P_k^\times, P^\times)$
 is a region for $\PP\cup\{\Pbold\}$.
\end{lemma}

\begin{proof} If $\Pbold$ commutes with $\Pbold_i$ for all $i$ either choice of sides will do.   So suppose $\Pbold$ does not commute with $\Pbold_i$ for some $i$.  By Lemma \ref{commute}, if $P_i^\times$ is not contained in either side of $\Pbold$, then either choice of side  works for this pair.  If $P_i^\times$ is contained in one side of $\Pbold$, we must choose $P^\times$ to be that side.  Thus, to prove the lemma, we must show that if $P_i^\times$ and $P_j^\times$ are each contained in a side of $\Pbold$, then they are contained in the same side.

To see this, suppose that $P_i^\times \subset P$ and $P_j^\times \subset P^*$, so  $P_i^\times \cap P_j^\times =\emptyset$.  By assumption,
the choice of sides for $\Pbold_i$ and $\Pbold_j$ defined a region, so they must commute.  That is, ${\rm max}(\Pbold_i)$ and ${\rm max}(\Pbold_j)$ are adjacent in $\G$.  Let $v \in {\rm max}(\Pbold^\times_i)$.  If $v\in  {\rm single}(\Pbold)$, then ${\rm max}(\Pbold_j)\subseteq {\rm lk}(v)\subseteq {\rm lk}(\Pbold)$,  contradicting the assumption that $\Pbold$ does not commute with $\Pbold_j$.    Thus, the elements of  ${\rm max}(\Pbold_i)$ appear as doubles in $P$.  Likewise, elements of  ${\rm max}(\Pbold_j)$ appear as doubles in $P^*$.  
 But since ${\rm max}(\Pbold_i), {\rm max}(\Pbold_j)$ are adjacent, they lie in the same component of $\G \setminus {\rm st}(v), v \in {\rm max}(\Pbold)$,
hence they must appear on the same side of $\Pbold$.  
\end{proof}

We continue building $\SGPP$ by attaching edges to $\mathbb E_\PP$ for each element of $V$.  
We need the  following lemmas  in order to explain how this is done.

Associated to a region $R=(P_1^\times,\dots, P_k^\times)$ is a subset of $V^\pm$ defined by  
 $$I(R) = \Pbar^{\times}_1 \cap \dots \cap \Pbar^{\times}_k $$
where $\Pbar^{\times}_i = P^{\times}_i \cup {\rm lk}(P_i)$.  As we will see below, the elements of $I(R)$ will correspond to the directed edges attached at the vertex $x_R$

Note that if switching sides of $\Pbold_i$ and leaving all other $P_j^{\times}$ unchanged gives a valid region $R_i$,  then there is an edge in $\mathbb E_\PP$ labelled $e_{\Pbold_i}$ from $x_R$ to $x_{R_i}$.   

\begin{lemma}\label{edges}
 Let $\PP= \{\Pbold_1, \dots ,\Pbold_k\}$ be a set of pairwise compatible \GW partitions.  Then the regions of $\PP$ satisfy the following.
 \begin{enumerate}
\item  For every element $v$  in $V^\pm$, there exists a region $R$ with $v \in I(R)$. 
\item If $I(R)$ contains $v$, then switching sides of every $\Pbold_i$ for which $v$ is a singleton gives a region $R_v$ such that $I(R_v)$ contains $v^{-1}$.  Moreover, if $w \in I(R)$ commutes with $v$, then $I(R_v)$ also contains $w$.  
\end{enumerate}
\end{lemma}

\begin{proof} 
We proceed by induction on $k$.  For $k=1$ this is obvious.  Assume $k>1$.  If $v \in {\rm lk}(P_i)$ for all $i$, then for \emph{every} region $R$,  $v \in I(R)$ and $R_v=R$, so we need only show that $\PP$ has at least one region, which follows from Lemma~\ref{regions}.  

Now suppose that 
$v$ is not in the link of at least one of the partitions.  Say $v \notin {\rm lk}(P_i)$ for $1\leq i < r $.  For these partitions, choose $P_i^\times$ to be the side of $\Pbold_i$ containing $v$.     For $r \geq i \leq k$, choose any collection of sides $P_i^\times$ such that $(P_r^\times, \dots ,P_k^\times)$ is a region for the partition $(\Pbold_{r}, \dots ,\Pbold_k)$.  (Such a choice exists by induction.)  We claim that $R=(P_1^\times, \dots ,P_k^\times)$ is a valid region for $\PP$.  For this, we must verify that the chosen sides for any two non-commuting partitions $\Pbold_i, \Pbold_j$ intersect. This is true by definition for $i,j \geq r$.  For $i,j < r$, $v \in P_i^\times \cap P_j^\times$ so this intersection is non-empty, and for $i < r \geq  j$,  $v \in P_i^\times \cap {\rm lk}(P_j)$ so by Lemma \ref{commute}, $\Pbold_i,\Pbold_j$ commute.  Thus $R$ is a region.  

If $R_v$  is obtained from $R$ by switching the sides of those $\Pbold_i$ for which $v$ is a singleton,
then the same argument, with $v$ replaced by $v^{-1}$, shows that $R_v$ is a region.  Moreover, if $v$ is a singleton in $\Pbold_i$, then 
${\rm lk}(v) \subseteq {\rm lk}(P_i)$, hence $I(R) \cap {\rm lk}(v)^\pm = I(R_v) \cap {\rm lk}(v)^\pm$.  
\end{proof}

 We can now explain how to attach edges to $\mathbb E_\PP$.  For each $v$  in $I(R)$ we attach an edge labelled $e_v$ joining $x_R$ and $x_{R_v},$ where $R_v$ is obtained as in the lemma.  (Note that we may have $R=R_v$ if $v$ is not a singleton in any $\Pbold_i$, in which case we are attaching a loop.)  The edge $e_v$ is oriented from the region containing $v^{-1}$ to the region containing $v$.  The 1-skeleton of the resulting complex, including edges labeled $e_v$ and $e_{\Pbold_i}$,  will be the entire 1-skeleton of $\SGPP$, so we denote it $(\SGPP)^{(1)}$

Note that while a given label occurs at most once at each vertex, it does not determine a unique edge in $(\SGPP)^{(1)}$.
For example, an edge labelled $e_v$ will occur at every  vertex $x_R$ with  $v \in I(R)$.  Indeed, once the higher dimensional cells are added, we will see that two edges have the same label if and only if they determine the same hyperplane in $\SGPP$.  
 
To complete the construction of $\SGPP$ we need to add higher-dimensional cubes which capture the commutation relations in $\AG$.   Define two edges to have {\em commuting labels} if their labels are one of the following.
\begin{enumerate}
\item  $e_v, e_w$ with $v,w$ distinct, commuting elements of $V,$
\item  $e_v, e_{\Pbold_i}$ with $v \in {\rm lk}(\Pbold_i)$,
\item  $e_{\Pbold_i}, e_{\Pbold_j}$  with $\Pbold_i, \Pbold_j$ distinct, commuting partitions.
\end{enumerate} 
 
\begin{lemma}\label{square}  Let $e_a$ and $e_b$ be edges at a vertex $x_R$ with commuting labels.  Then $e_a,e_b$ belong to a 4-cycle in $(\SGPP)^{(1)}$ with opposite edges having the same label.
\end{lemma}

\begin{proof}
If the labels are both of the form $e_{\Pbold_i}$ then they span a square in $\mathbb E_\PP$, and we are done. 

If the labels are   $e_v$ and $e_w$, they terminate at $x_{R_v}$ and $x_{R_w}$ respectively.  It follows from Lemma \ref{edges}(2) that there is an edge labelled $e_v$ emanating from $x_{R_w}$ and an edge labelled $e_w$ emanating from $x_{R_v}$.  These form a square with the vertex opposite $x_R$ corresponding to $(R_v)_w = (R_w)_v$.

A similar argument applies for edges labelled $e_v$ and $e_{\Pbold_i}$.  Since $v \in {\rm lk}(\Pbold_i)$, switching the side of $\Pbold_i$ does not affect $v$.  So if $R^\prime$ is the result of this switch, then there is an edge labelled $e_v$ emanating from $x_{R^\prime}$. The other end of this edge corresponds to a region $R^\prime_v$ which differs from $R_v$ only on $\Pbold_i$.  Thus,  the vertices $R_v$ and $R^\prime_v$ are also connected by an edge labelled $e_{\Pbold_i}$, completing the square.
\end{proof}

\begin{corollary}  If a collection of edges $e_1, \dots e_m$ emanating from a vertex $x_R$ have pairwise commuting labels, then these edges form a corner  of the 1-skeleton of an $m$-cube  in $(\SGPP)^{(1)}$, such that parallel edges have the same labels.  
\end{corollary}

\begin{proof}  This follows from Lemma~\ref{square} by induction on $m$.
\end{proof}

It follows that we can glue an $m$-cube into $(\SGPP)^{(1)}$ whenever we have a set of $m$  edges at a vertex with commuting labels.   The resulting cube complex is  $\SGPP$.   Note that the subcomplex of $\SGPP$ spanned by the edges labelled $e_{\Pbold_i}$ is precisely the complex $\mathbb E_\PP$ that we started with.   By construction,  the link of the vertex $x_R$ has an $(m-1)$--simplex for each set of $m$ mutually commuting edge labels emanating from $x_R$, so by Gromov's link condition $\SGPP$ is locally CAT(0). 

\begin{definition}  A {\em $\G$-complex} is any  cube complex which is isomorphic to $\SGPP$ for some compatible set $\PP$ of $\G$-Whitehead partitions, including $\SG=\SG^{\emptyset}$.  A $\G$-complex isomorphic to $\SG$ is called simply a {\em Salvetti.}  
\end{definition}

\begin{theorem}\label{thm:blowup} Let $\PP$ be a compatible set of \GW partitions.  Then the $\G$-complex   $\SGPP$  has the following properties.
\begin{enumerate}
\item  Any two vertices, $x_R, x_{R'}$ are connected by a path with labels in the set 
$$\{e_{\Pbold_i} \mid \textrm{$R$ and $R'$ contain opposite sides of $\Pbold_i$ }\}. $$
\item Any two edges with the same label are dual to the same hyperplane.
\item $\SGPP$ is connected and locally CAT(0).
\end{enumerate}
\end{theorem}

\begin{proof}  (1)  Let $R$ and $R^\prime$ be two regions of $\PP$ and (reordering if necessary) suppose they differ in the choice of sides of $\Pbold_1, \dots ,\Pbold_l$ but agree on the remaining partitions. We will show by induction on $l$ that they are connected by a path with  labels in $\{e_{\Pbold_i}\}_{i=1, \dots , l}$.
 For $l=1$, two regions which differ on only one partition $\Pbold$ are, by construction, connected by an edge labelled $e_{\Pbold}$. 

Suppose $l >1$.    For simplicity, write $R=(P_1, \dots ,P_k)$.  (It makes no difference which side we call $P_i$ and which side we call $P_i^*$).  Choose $i \leq l$ such that $P_i$ is minimal among the sets $P_1, \dots ,P_l$.  That is, $P_i$ does not contain any of the other sets in this collection. Then it follows from Lemma \ref{commute}  that switching $P_i^*$  is also compatible with (i.e. commutes or intersects)  the remaining $P_j, j \leq l$.  For $j> l $, $P_j$ must also be compatible with  $P_i^*$ since they both appear in the region $R^\prime$.  Thus, setting $R_i = (P_1, \dots ,P_i^*, \dots P_k)$ we obtain a valid region which differs from $R^\prime$ in $l-1$ places.  The vertices corresponding to $R$ and $R_i$ are connected by an edge labelled $e_{\Pbold_i}$ and by induction, the vertices corresponding to $R_i$ and $R^\prime$ are connected with  labels in $\{e_{\Pbold_j}\}_{1\leq j \leq l, j \neq i}$.

(2)  We consider two cases. Suppose $R$ and $R^\prime$ are regions with $v \in I(R) \cap I(R^\prime)$, so there is an edge labelled $e_v$ emanating from  both vertices $x_R$ and $x_{R^\prime}$.  Then for any partition $\Pbold_i$ with $v \notin {\rm lk}(\Pbold_i)$ both $R$ and $R^\prime$ must contain the (unique) side of $\Pbold_i$ containing $v$.  That is, $R$ and $R^\prime$ differ only on partitions containing $v$ in their link.  By part (1), it follows that $x_R$ and $x_{R^\prime}$ are connected by a path labeled by $e_{\Pbold_i}$ such that $e_{\Pbold_i}$ commutes with $e_v$ and hence these two edges span a 2-cube. Proceeding along this path gives a sequence of such cubes joining the $e_v$ edges at $x_R$ and $x_{R^\prime}$.  It follows that they are dual to the same hyperplane.

For two edges labelled by $e_{\Pbold_i}$  consider the four vertices contained in these edges.  Say the regions for these vertices are $R, R_i$ and $R^\prime,R^\prime_i$.  If $R$ and $R^\prime$ differ on some $P_j$ then all possible combinations of $P_i,P_i^*$ with $P_j,P_j^*$ occur in these four regions.  But this is possible only if $\Pbold_i$ and $\Pbold_j$ are commuting partitions.  Thus, arguing as above, we can connect $R$ to $R^\prime$ with a path labelled by $e_{\Pbold_j}$'s which commute with $e_{\Pbold_i}$ and conclude that there is a sequence of cubes between the $e_{\Pbold_i}$-edges at $x_R$ and $x_{R^\prime}$.

(3)  It follows from (1) that $\SGPP$  is connected, and it was observed above that it is locally CAT(0) by construction.   

\end{proof}


\section{Collapsing along hyperplanes}

In the case of a single blow-up $\SGP$, we observed in Lemma \ref{blowupcollapse} that for any element  $m\in {\rm max}(P)$ there is a subcomplex containing the edge $e_m$ which can be collapsed to give a Salvetti, i.e. a cube complex isomorphic to $\SG$.  Furthermore the map on $\SG$ obtained by blowing-up followed by this collapse corresponds to a Whitehead automorphism.  
In this section we identify all subcomplexes of general $\G$-complexes which can be collapsed to give a Salvetti.    

Let $X$ be a  non-positively curved cubical complex and $H$ a hyperplane of $X$.   If $e$ is an edge of $X$ then  $e$ and $H$ are said to be \emph{dual} if $e$ intersects $H$. 
The \emph{carrier} $\kappa(H)$ of $H$ is the subcomplex of $X$ formed by the closures of the cubes of all dimensions that intersect $H$.  

\begin{definition} Let $H$ be a hyperplane of $X$.  We say   $H$ is a \emph{carrier retract} if $\kappa(H)$ is isomorphic to $H \times [0,1]$; 
in particular $H$ is embedded in $X$ and there are no identifications on the boundary of $\kappa(H)$.  If $H$ is a carrier retract, we define the \emph{collapse of $X$ along $H$} to be the cube complex formed by collapsing $\kappa(H)$ orthogonally onto $H$.  Denote the resulting complex by $X_H$, and note that there is a canonical projection $X \twoheadrightarrow X_H$.
\end{definition}
\begin{example}  In the blowup $\SGP,$ there is one hyperplane dual to $e_\Pbold$ and one dual to $e_v$ for each $v\in V$.   The hyperplane dual to $e_\Pbold$ is isomorphic to $\SG_{{\rm lk}(\Pbold)}$, and   for every $v\in V$  the hyperplane dual to $e_v$  is  isomorphic to $\SG_{{\rm lk}(v)}$.    The hyperplane dual to $e_\Pbold$ is a carrier retract. The hyperplane dual to $e_v$ is a carrier retract if and only if $v\in {\rm single}(\Pbold)$.
\end{example}  

\begin{definition} Let $\mathcal H= \{H_1, H_2, \dots, H_k\}$ be a set of hyperplanes.  We say that $\mathcal H$  \emph{has compatible carriers} if  each $H_i$ is a carrier retract and any loop in $X$ consisting of edges dual to the $H_i$'s is null homotopic.  Given such a set, define the \emph{collapse of $X$ along $\mathcal H$} to be the complex $X_{\mathcal H}$ obtained by collapsing each cube $C$ in $\bigcup \kappa(H_i)$ to the intersection of the mid-planes of $C$ lying in some $H_i$. 
 \end{definition}

 The proof of the following lemma is an easy exercise.
 
 \begin{lemma}  Let  $\mathcal H= \{H_1, H_2, \dots ,H_k\}$ be set of hyperplanes in $X$ with compatible carriers.  Let $Y=X_{H_1}$ and for $i>1$, let $\overline H_i$ denote the image of $H_i$ in $Y$.  Then
 \begin{enumerate}
 \item $Y$ is a non-positively curved cube complex.
 \item  $\overline{\mathcal H}= \{\overline H_2, \dots ,\overline H_k\}$ has compatible carriers in $Y$.
 \item $Y_{\overline{\mathcal H}} = X_{\mathcal H}$.
 \end{enumerate}
 \end{lemma}

In particular, it follows that collapsing along the $H_i$'s one at a time in any order results in the same space $X_{\mathcal H}$.  

\begin{lemma}\label{minimal paths}  Let  $\mathcal H$ be   set of hyperplanes in $X$ with compatible carriers and let $c : X \to X_{\Hcal}$ be the projection map.  Then
\begin{enumerate}
\item $c$ is a homotopy equivalence, 
\item  distinct hyperplanes in $X$ not contained in $\Hcal$ map to distinct hyperplanes in $X_\Hcal$,
\item if $g$ is a conjugacy class in $\pi_1(X)$ and $p$ is a minimal length edge path in $X$ representing $g$, then $c(p)$ is a minimal length edge path in $X_\Hcal$ representing $g$.
\end{enumerate}
\end{lemma}

\begin{proof}  By the previous lemma, it suffices to consider the case where $\Hcal$ consists of a single hyperplane $H$ which is a carrier retract.  The first statement is clear:  the homotopy equivalence between $\kappa(H)$ and $H$ extends to a homotopy equivalence between $X$ and $X_H$ since $\kappa(H)$ is a strong deformation retract of the open neighborhood  consisting of points at distance $<\frac{1}{2}$ from $\kappa(H)$. 

  For the second statement, recall that a hyperplane can be identified with an equivalence class of edges.  Two edges in $\kappa(H)$ that become identified under the collapse $c$ are parallel edges in some cube, hence they are already equivalent.  

For the third statement, let $g$ be a conjugacy class in $\pi_1(X)$ and let $p$ be an edge path in $X$ representing $g$.  Lift $p$ to a path $\tilde p$ in the universal cover $\widetilde X$ and let $\tilde p^\infty$ denote the union of the $g^k$-translates of $\tilde p$, for $k\in\mathbb Z$.  Since $\widetilde X$ is a CAT(0) cube complex, it follows from \cite{Vij} that $p$ is minimal if and only if $\tilde p^\infty$ crosses no hyperplane of $\widetilde X$  more than once. The analogous statement holds for $c(p)$.  

The universal cover of $X_\Hcal$ is obtained from $\widetilde X$ by collapsing along  all  hyperplanes $\widetilde \Hcal$ in the inverse image of $\Hcal$.  Let $\tilde c : \widetilde X \to \widetilde X_{\widetilde \Hcal}$ be the lift of $c$.  Set $q=c(p)$,  and define $\tilde q^\infty$ as above.  Then $\tilde c(\tilde p^\infty) = \tilde q^\infty$, so
by part (2) of the lemma, if $\tilde p^\infty$ crosses each hyperplane at most once, the same holds for 
$ \tilde q^\infty$.  Statement (3) follows.
\end{proof}

We are now ready to apply these observations to the hyperplanes of a blowup $\SGPP$. 
\begin{theorem}  Let $\PP=\{\Pbold_1,\ldots,\Pbold_k\}$ be a compatible set of \GW partitions for $\Gamma$.  Then the  set of hyperplanes $\mathcal H = \{H_1, \dots ,H_k\}$ dual to the edges $e_{\Pbold_i}$ has compatible carriers.  For any subset $J \subseteq \{1, 2, \dots , k\}$, collapsing along the hyperplanes $\{H_i\}_{i \in J}$  gives a complex naturally isomorphic to the blowup of $\SG$ by the partitions $\{\Pbold_i\}_{i \notin J}$.  In particular, collapsing along all of $\mathcal H$ gives a complex isomorphic to the Salvetti complex $\SG$.
\end{theorem}

\begin{proof}    It is easy to check that the carrier of each $H_i$ is isomorphic to $H_i\times [0,1]$.  Thus to show that $\mathcal H$ has compatible carriers, we need to show that any loop $\gamma$ formed by edges labelled $e_{\Pbold_i}$ is null homotopic. We will induct on the number of $\Pbold_i$'s appearing in this loop.  Say the loop involves only a single $\Pbold_i$.  Since no two edges with the same label occur at a vertex, the loop would have to consist of a single edge.  But an $e_{\Pbold_i}$-edge switches the side of $\Pbold_i$, so it cannot be a loop.  

Now suppose $\gamma$ involves more than one $\Pbold_i$.  Orient $\gamma$ and say the initial vertex is $x_R$ with $R=(P_1, \dots ,P_k)$.  Once $\gamma$ crosses an $e_{\Pbold_i}$-edge, all of the regions it encounters will have $P_i^*$ in the $i^{th}$ position until it crosses another $e_{\Pbold_i}$-edge.  Thus, the labels must occur in pairs. Let $\alpha$ be a segment of $\gamma$ joining two consecutive $e_{\Pbold_i}$-edges, so $\gamma$ decomposes as
$$ \gamma = \gamma_1\cdot e_1 \cdot \alpha\cdot e_2 \cdot \gamma_2. $$
where $e_1,e_2$ are edges labeled by $e_{\Pbold_i}$.
By Theorem~\ref{thm:blowup}(2), there exists a path $\beta$ in the carrier $\kappa(H_i)$ between the endpoints of $\alpha$.  Then $\alpha\beta^{-1}$ forms a loop not involving $\Pbold_i$, so by induction $\alpha$ is homotopic 
to $\beta$.  Moreover, the path $e_1 \beta e_2$ can be slid across the hyperplane $H_i$ to get a path $\beta'$ with no edges labelled $e_{\Pbold_i}$.
Thus, $\gamma$ is homotopic to $\gamma^\prime=\gamma_1\beta'\gamma_2$.  Repeating this process if necessary, we can get rid of all $e_{\Pbold_i}$-edges in $\gamma$ and apply induction to conclude that $\gamma$ is null-homotopic.

For the second statement, it suffices to consider the case where $J$ is a singleton, say $J=\{1\}$.  So consider the space obtained from $\SGPP$ by collapsing along $H_1$, the hyperplane dual to $e_{\Pbold_1}$.  This collapse has the effect of identifying two vertices whose label differ only in the  choice of side for $\Pbold_1$.  So letting
$\PP'=\{\Pbold_2, \dots ,\Pbold_k\}$,  we can map vertices of the quotient space injectively to vertices of $S^{\PP^\prime}$  by forgetting $\Pbold_1$.  By Lemma \ref{regions}, this map is also surjective.  The construction of $S^{\PP^\prime}$  depends only on the vertex labels, so it is now easy to verify that this bijection extends to an isomorphism of complexes.  
\end{proof}

We will call any collapse along hyperplanes dual to $e_{\Pbold_i}$-edges, $i \in J \subseteq \{1, \dots ,k\},$ a \emph{canonical collapse} of $\SGPP$.
In particular, taking $J=\{1, \dots ,k\}$, we have a canonical collapse from $\SGPP$ down to the Salvetti complex $\SG$.  However, one can obtain a Salvetti complex by collapsing along many  other sets of hyperplanes $\mathcal H$  in  $\SGPP$. 

\begin{example}\label{Whitehead auto}  
Let $\Pbold =\{P,P^*, {\rm lk}(P)\}$ be a single \GW partition and $\SGP$ the associated $\G$-complex.  As we saw in Lemma~\ref{blowupcollapse}, collapsing along the hyperplane dual to the edge labelled $e_v$   for any $v \in {\rm max}(P)$ gives a complex isomorphic to $\SG$. 
\end{example}

\begin{example}\label{pair} 
Let $\Pbold$ and $\Qbold$ be   compatible $\Gamma$-partitions with ${\rm lk}(\Pbold)={\rm lk}(\Qbold)=L.$    
In particular $\Pbold$ and $\Qbold$ do not commute, so $\SG^{\{\Pbold,\Qbold\}}$ has one edge labeled $e_\Pbold$ and one edge labeled $e_\Qbold$.  Let $\Theta$ denote the graph formed by $e_\Pbold$, $e_\Qbold$ and all $e_v$ with ${\rm lk}(v)=L$.   The hyperplane dual to each  edge in $\Theta$ is isomorphic to the Salvetti complex $\SG_L$.  Thus the subcomplex  spanned by the carriers of all of these hyperplanes decomposes as a product $\Theta \times \SG_L$.   Now take any maximal tree $T$ in $\Theta$ and let  $\mathcal H$  be the set of hyperplanes dual to the edges in $T$.  Then collapsing $\SG^{\{\Pbold,\Qbold\}}$ along $\mathcal H$ reduces $\Theta$ to a rose and reduces $\SG^{\{\Pbold,\Qbold\}}$ to a complex isomorphic to $\SG$.  
\end{example}

\begin{example}\label{samelink} 
Example ~\ref{pair} generalizes to any set $\PP$  of compatible $\Gamma$-partitions which all have  the same link $L$.  Since no two elements of $\PP$ commute, the 1-skeleton of $\SG^\PP$ has exactly one edge labeled $e_\Pbold$ for each $\Pbold\in\PP$ and exactly one labeled $e_v$ for each $v$ which is not in $L$.  For every $v$ the hyperplane dual to $e_v$ is isomorphic to $\SG_{{\rm lk}(v)}$, and the hyperplane dual to each $e_\Pbold$ is isomorphic to $\SG_L$. If   $[m]$ denotes the set of all vertices $v\in\G$ with ${\rm lk}(v)=L$, then the union of the carriers of the hyperplanes dual to the edges $e_{\Pbold}$ for $\Pbold\in \PP$ and $e_v$ for $v\in [m]$ decomposes as a product $\Theta \times \SG_L$.  We will call $\Theta$ the \emph{base graph} of $\PP$. Collapsing $\SGPP$ along any set of hyperplanes dual to a maximal tree in the base graph reduces $\SGPP$ to a complex isomorphic to the Salvetti complex $\SG$ for $\Gamma$.

Note that a hyperplane $H$ of $\SG^\PP$  is a carrier retract if and only if the dual edge $e$ to $H$ is not a loop, i.e. if and only if $e=e_\Pbold$ for some $\Pbold$ or $e=e_v$ for $v$ a singleton in some $\Pbold$.   A set $\mathcal K=\{H_e\}$ of hyperplanes has compatible carriers if and only if the dual edges form a forest in the 1-skeleton.  
\end{example}

Now  let $\PP$ be any set of compatible \GW partitions. Subdivide $\PP$ into subsets $\PP= \PP_1 \cup \dots \cup \PP_s$ where each $\PP_i$ is a maximal collection of $\Pbold_j$ having the same link, $L_i$.  (For example, if $\AG$ is a free group then all of the links are empty, hence $s=1$ and $\PP=\PP_1$.)  Consider the blow-ups $\SG^{\PP_i}$.  By the discussion above, each of these contains a subcomplex of the form $\Theta_i \times \SG_{L_i}$ where $\Theta_i$ is the base graph of $\PP_i$.  The edges of $\Theta_i$ are labelled by $e_\Pbold$ and $e_v$ with $\Pbold \in \PP_i$ and ${\rm lk}(v)=L_i$.  In particular,  for $i \neq j$, the labels on the edges of $\Theta_i$ and $\Theta_j$ are disjoint.

\begin{definition}   Let $\PP$ be a set of compatible \GW partitions.  If all of the partitions in $\PP$ have the same link, call a set of hyperplanes $\Hcal$ in $\SGPP$  \emph{tree-like} if the edges dual to $\Hcal$ form a maximal tree in the base graph $\Theta$.  More generally, if  $\PP= \PP_1 \cup \dots \cup \PP_s$ is the decomposition into partitions with the same link,  call  $\mathcal H$  \emph{tree-like} if $\Hcal = \Hcal_1 \cup \dots \cup \Hcal_s$ where
 $\Hcal_i$ is the lift of a tree-like set of hyperplanes in $\SG^{\PP_i}$, that is, the edges dual to $\Hcal_i$ form a maximal tree in $\Theta_i$.
\end{definition}

For example, the set of hyperplanes dual to the edges labelled $e_{\Pbold_i}, 1 \leq i \leq k$ is always tree-like.

\begin{theorem}\label{treelike}  Let $\Hcal$ be a set of  hyperplanes in $\SGPP$.  Then  $\mathcal H$ is tree-like  if and only if it has compatible carriers and the collapse of $\SGPP$ along $\mathcal H$ is isomorphic to the Salvetti complex $\SG$.
\end{theorem}

\begin{proof}  Suppose $\Hcal=\Hcal_1\cup\ldots\cup\Hcal_s$ is tree-like.   We proceed by induction on $s$; the case $s=1$ was done in  Example \ref{samelink}.
For $s>1$, reordering if necessary we may assume that the link $L_1$ is maximal among the  $L_i$'s.  Suppose $v$ is a vertex of
$\G$ with ${\rm lk}(v)=L_1$. If $v$ is a singleton in some $\Pbold \in \PP$, then $v \leq {\rm max}(\Pbold)$ implies that $L_1 \subseteq {\rm lk}(\Pbold)$,  so by our maximality assumption $L_1 = {\rm lk}(\Pbold)$ and hence  $\Pbold \in \PP_1$.  It follows that edges labelled  $e_v$ connect vertices of $\SGPP$ which differ only on partitions  in  $\PP_1$, and hence the graph $\Theta_1$ lifts isomorphically to a graph with the same labels in $\SGPP$.  The carriers of the hyperplanes dual to this graph in  $\SGPP$ span a subcomplex $Y= H \times \Theta_1\iso \SG_{L_1}\times \Theta_1$. 

Since $\Hcal$ is tree-like, $\Hcal_1$ consists of hyperplanes in $Y$ dual to some maximal tree in $\Theta_1$. In particular, $\Hcal_1$ has compatible carriers and the resulting collapse reduces the subcomplex $Y$  to the product of $H$ with a rose and leaves everything else unchanged.  The resulting complex is thus isomorphic to 
$\SG^{\PP'}$ where $\PP'=\PP \setminus \PP_1$.  Let $\rho : \SGPP \to S^{\PP'}$ be the collapsing map followed by this isomorphism.

The image of the hyperplanes $\Hcal' = \Hcal \setminus \Hcal_1$ in $\SG^{\PP'}$ is tree-like (since the canonical projection from $\SGPP$ to $\SG^{\PP_i}$ factors through $\SG^{\PP'}$), so by induction, it has compatible carriers and the resulting space is isomorphic to the Salvetti complex $\SG$.  It now follows that the original set $\Hcal$ has compatible carriers in $\SGPP$ since if $p$ is a loop of edges dual to $\Hcal$, then its image in $S^{\PP'}$ is a loop dual to $\Hcal'$.  This loop must be null-homotopic hence the same holds for $p$. 

Conversely, suppose $\Hcal$ is a set of hyperplanes in $\SGPP$ which has compatible carriers and collapses  $\SGPP$ down to $\SG$.  We again proceed by induction on $s$.  The case  $s=1$ is discussed in Example~\ref{samelink}, where it is observed that since $\Hcal$ has compatible carriers, the dual edges $\{e_i\}$ form a forest in the 1-skeleton of $\SG^\PP$.  Since collapsing along $\Hcal$ reduces $\SG^\PP$ to a complex isomorphic to $\SG$, which has only one vertex, these edges must form a maximal tree $T$ in the 1-skeleton. Since edges of $T$ are not loops, they correspond to the $\Pbold_i$ or to singletons in the $\Pbold_i$, so the hyperplane dual to each edge in $T$ is isomorphic to a subcomplex of $\SG_L$.  An edge of $T$ cannot correspond to a non-maximal singleton $v$,  since then ${\rm lk}(v)$ would be a proper subcomplex of $L$,  the carrier of $\Hcal$ would have fewer cubes than the carrier of $\{H_{\Pbold}\}_{\Pbold\in\PP}$, and collapsing along $\Hcal$ and along $\{H_{\Pbold}\}_{\Pbold\in\PP}$ would not result in isomorphic complexes.  Therefore $T$ is a maximal tree in the base graph, i.e. $\Hcal$ is treelike.

Now suppose $s>1$.  As above, assume that the link $L_1$ is maximal, so the graph $\Theta_1$ may be viewed as a subcomplex of $\SGPP$.  
Let $\Hcal_1$ be the set of hyperplanes in $\Hcal$ dual to some edge of $\Theta_1$.  We claim that these edges form a maximal tree in $\Theta_1$. 
Let $Z = (\SGPP)_\Hcal$ be the collapse of $\SGPP$ along $\mathcal H$ and 
$c_{\mathcal H} : \SGPP \to Z$ the collapsing map.  Let
$c : \SGPP \to \SG$ be the canonical collapse. 
By assumption, $Z$ is isomorphic to the Salvetti complex $\SG$,  so the image of $\Theta_1$ under  both $c$ and  $c_{\mathcal H}$ is a rose.    The former generates a free subgroup of the fundamental group $\AG$ (namely the subgroup generated by the vertices of $\G$ with link equal to $L_1$).  
Since both collapsing maps are homotopy equivalences, the same must be true of the latter.  It follows that the edges dual to $\Hcal_1$ must form a maximal tree in $\Theta_1$. 

Now set  $X=(\SGPP)_{\Hcal_1}$, and $\PP'=\PP \setminus \PP_1$.  Then the images of $\Theta_1$ in $S^{\PP'}$ and in $X$ are isomorphic roses and this isomorphism extends to an isomorphism of the whole complex $S^{\PP'} \cong X$.  Collapsing $X$ along the image  $\Hcal' $ of  $\Hcal \setminus \Hcal_1$ gives $Z$, hence by induction, it is tree-like (viewed as hyperplanes in $S^{\PP'}$).  We conclude that  that the original set of hyperplanes $\Hcal$ was tree-like in $\SGPP$.
\end{proof}

\begin{theorem}\label{swap}  Let $\Hcal$ and $\Kcal$ be two  tree-like sets of hyperplanes in $\SGPP$.  Given any $K\in \Kcal$, there exists  $H \in {\Hcal}$ such that the set of hyperplanes obtained from $\Hcal$ by replacing $H$ by $K$ is again tree-like.
\end{theorem}

\begin{proof}  Note that the label dual to a hyperplane in $\Hcal$ (or $\Kcal$) appears in one and only one of the graphs $\Theta_i$ since each $\Theta_i$ corresponds to a different  link.
  Say $K \in \Kcal$ is dual to an edge $e$ in $\Theta_i$. Let $T_i$ be the maximal tree in $\Theta_i$  formed by edges dual to $\Hcal$.  If $e$ lies in $T_i$, then $K$ also lies in $\Hcal$ and we can take $H=K$. If not,  let $e'$ be an edge in $T_i$ on the minimal path between the two vertices of $e$.  Then replacing $e'$ by $e$ gives another maximal tree $T_i'$ in $\Theta_i$ and has no effect on the remaining $\Theta_j$.  Thus, replacing the hyperplane $H$ dual to $e'$ by $K$ gives another tree-like set.
\end{proof}

\begin{corollary}\label{independent}
Let $\PP=\{\Pbold_1,\ldots,\Pbold_k\}$ be a compatible collection of \GW partitions, let $H_i$ be the hyperplane in $\SGPP$ dual to the edges $e_{\Pbold_i}$, and let   $\Kcal$ be another a treelike set of hyperplanes in $\SGPP$.  Then the automorphism of $\AG$ induced by $\SG\leftarrow \SGPP\to (\SGPP)_{\Kcal}\iso\SG$ is an element of $\UAG$.   
\end{corollary}

\begin{proof} By Theorem~\ref{swap}, we can order the elements of $\Hcal$ as $\{H_1,\ldots,H_k\}$ so that for each $i=0,\ldots,k$ the set $\Hcal_{i}=\{H_1,\ldots,H_i,K_{i+1},\ldots, K_k\}$ is treelike.   If we now set $$\widehat\Hcal_i=\{H_1,\ldots,H_{i-1},K_{i+1},\ldots,K_k\} \hbox{ for $i=1,\ldots,k$}$$ then the blowup-collapse $\SG\leftarrow \SGPP\to (\SGPP)_\Kcal$ factors into the sequence of single blowup-collapses
$$\SG=(\SGPP)_{\Hcal_k}\leftarrow (\SGPP)_{\widehat\Hcal_k}\to (\SGPP)_{\Hcal_{k-1}}\leftarrow\ldots    
       (\SGPP)_{\Hcal_{k-1}}\leftarrow (\SGPP)_{\widehat\Hcal_i}\to (\SGPP)_{\Hcal_{i}}\ldots \to $$ $$\ldots \to (\SGPP)_{\Hcal_{1}}\leftarrow (\SGPP)_{\widehat\Hcal_{1}}\to (\SGPP)_{\Hcal_{0}}=(\SGPP)_\Kcal$$
The statement of the corollary now follows from Lemmas~\ref{induced} and~\ref{blowupcollapse}.
\end{proof}

\subsection{Construction of \texorpdfstring{$\KG$}{KG}}

We are now ready to  define the simplicial complex $\KG$ as  the geometric realization of a partially ordered set  of  $\G$-complexes.  
 
\begin{definition} A \emph{marked $\G$-complex} $\sigma$ is a pair,  $\sigma=(X,\alpha)$ where 
\begin{enumerate}
\item $X$ is a cube complex isomorphic to $\SGPP$ for some (possibly empty) compatible set $\PP$ of \GW partitions.  
\item  $\alpha \colon X\to \SG$ is a homotopy equivalence and the composition $\SG  \buildrel{c_\PP^{-1}}\over\to \SGPP \iso X\buildrel{\alpha}\over\to\SG$ induces an element of $\UAG$. 
\end{enumerate}
Two marked $\G$-complexes $\sigma=(X,\alpha)$ and $\sigma'=(X^\prime,\alpha^\prime)$ are \emph{equivalent} if there is an isomorphism of cube complexes $h\colon X\to X^\prime$ with $\alpha^\prime\circ h\simeq \alpha$. 
If $X$ is isomorphic to $\SG$, the equivalence class of $(X,\alpha)$ is called  a \emph{\rose}.
\end{definition}
Note that the second condition in the definition of a marked $\G$-complex is independent of the choice of isomorphism $X\iso \SGPP$ or collapse $c_\PP$ by Corollary~\ref{independent}.

\begin{examples}\label{isometries} (1) For a \GW pair $(P,m)$ we observed in the discussion preceding Lemma \ref{blowupcollapse}, that the collapsing maps $c_P$  and $c_m$ on $\SGP$ differ by the isomorphism that interchanges the hyperplanes dual to $e_P$ and $e_m$.  It follows that $(\SGP,c_P)  \sim (\SGP,c_m)$.

(2) If $\varphi \in \Out$ is a product of symmetries and inversions, then it can be represented by an isomorphism  $\hat\varphi \colon \SG \to \SG$, hence $(\SG, id) \sim (\SG, \hat\varphi). $
\end{examples}

We now define a partial ordering on the set of marked $\G$-complexes.
 If $\sigma=(X,\alpha)$,  $\Hcal$ is a set of hyperplanes of $X$ contained in some tree-like set, and $c : X \to X_\Hcal$ is the collapsing map, we denote by $\sigma_{\Hcal}$ the marked $\G$-complex $(X_\Hcal, \alpha \circ c^{-1})$.  For two marked $\G$-complexes $\sigma,\sigma'$, define
$$\textrm{$\sigma' < \sigma$  if $\sigma'=\sigma_\Hcal$ for some $\Hcal$.}$$

\begin{definition} The {\em $\G$-spine}  $\KG$ is the simplicial complex associated to the partially ordered set of equivalence classes of marked $\G$-complexes.  
\end{definition}

We can identify $\Out$ with the group of homotopy classes of maps $\SG \to \SG$.  Using this identification, we define a left action of $\UAG$ on $\KG$  by $\varphi \cdot (X, \alpha) = (X, \varphi \circ \alpha)$. 
\begin{proposition}  The action of $\UAG$ on $\KG$ is proper.
\end{proposition}  

\begin{proof}  Since each marked blowup $(X,\alpha)$ can be collapsed to finitely many marked Salvettis, it suffices to prove that the stabilizer of some (hence any) marked Salvetti is finite.  This is true for $(\SG, id)$ since any isomorphism $\SG \to \SG$ takes the one skeleton to the one skeleton, hence induces a permutation on $V^\pm$.  Thus $(\SG, id) \sim (\SG, \alpha)$ if and only if $\alpha$ lies in the (finite) group generated by graph symmetries and inversions.
\end{proof}

 Let $(P,m)$ be a \GW pair and $\alpha$ the corresponding Whitehead automorphism.  By Lemma~\ref{blowupcollapse}, $\alpha$ is realized by the blow-up-collapse, 
 $\alpha=c_m \circ c_\Pbold^{-1}\colon \SG\to \SGP\to  \SG$. 
   If we start at the Salvetti $(\SG, id)$ this gives a path in $\KG$ which ends at $(\SG,\alpha)$:  
\begin{equation}
(\SG, id) < (\SGP, c_P ) \sim  (\SGP, c_m ) = (\SGP, \alpha \circ c_P) > (\SG, \alpha). 
\end{equation}
More generally, for any $\varphi \in \UAG$, we can translate this path by $\varphi$ to obtain a path from 
$(\SG, \varphi)$ to $(\SG, \varphi\circ\alpha) $.
\begin{definition}  If $\sigma = (\SG, \varphi)$, we call the $\varphi$-translate of path (1) above the \emph{Whitehead move} at $\sigma$  associated to $(\Pbold,m)$, and write $\sigma^\Pbold_m= (\SG, \varphi\circ\alpha).$
\end{definition}

Using this terminology, Corollary~\ref{independent} can be restated in the following useful form:

 \begin{corollary}\label{factorization}  (Factorization Lemma)  Let $\sigma=(\SG,\alpha)$ be a \rose, $\PP=\{\Pbold_1,\ldots,\Pbold_k\}$ a compatible collection of \GW partitions,  $\sigma^\PP = (\SGPP, 
 c_{\PP} \circ \alpha)$ be the blow-up of $\sigma$ with respect to $\PP$, and $\Hcal$ a tree-like set of hyperplanes in $\SGPP$.  Then with a suitable ordering of the elements of $\Hcal$ there is a chain
 $\sigma=\sigma_0, \sigma_1,\ldots,\sigma_k=\sigma^\PP_\Hcal$ such that each $\sigma_i$ is connected to  $\sigma_{i-1}$ by a  Whitehead move.  
 \end{corollary}

Our goal is show that $\KG$ is contractible.  As a first step we have the following.  
\begin{proposition}\label{connected}  $\KG$ is connected.
\end{proposition}

\begin{proof} By definition, every vertex of $\KG$ lies in the star of some \rose.  It is straightforward to verify that the subgroup generated by \GW  automorphisms is normal in $\UAG$, hence any $\varphi \in \UAG$, can be factored   as a product $\varphi = \varphi_1 \circ \varphi_2$ where $\varphi_1$ is a product of symmetries and inversions  and $\varphi_2$ is a product of \GW automorphisms. It follows from Example~\ref{isometries} and Corollary~\ref{factorization} that $(\SG, id)=(\SG,\varphi_1)$ is connected by a path in $\KG$ to $(\SG,\varphi_1 \circ \varphi_2 )=(\SG, \varphi)$. 
\end{proof}


\section{Examples}\label{examples}

In this section we pause to compute a few examples of the complexes $\KG$.  These complexes quickly become very complicated and difficult to describe, but invariants such as dimension are relatively easy to compute.

In the next section we will show that  $\KG$ is contractible.  Since the action is proper, the dimension of $\KG$ gives an upper bound for the virtual cohomological dimension of $\UAG$.  We remark that the  kernel $I\AG$ of the natural map $\Out\to {\rm GL}(n,\Z)$ is a torsion-free subgroup of $\UAG$ so acts freely on $\KG$ and the dimension of $\KG$ is an upper bound for its cohomological dimension.  

Since there are only finitely many ways to partition a finite set, and maximal simplices of $\KG$ correspond to maximal sets of pairwise-compatible $\G$-partitions, $\KG$ is always finite-dimensional.  In any explicit example the dimension of $\KG$ can be computed precisely by finding a maximal set of pairwise-compatible $\G$-partitions.

\subsection{3-vertex graphs}
We first consider graphs with 3 vertices.

For any $n$, if $\G$ has $n$ vertices and no edges, then $\AG$ is the free group of rank $n$, the Salvetti complex is a rose with $n$ petals, and $\G$-complexes are graphs of rank $n$ without univalent or bivalent vertices and with no separating edges; in other words   $\KG$ is the spine  $K_n$ of (reduced) Outer space.  For $n=2$ this can be identified with the usual tree for ${\rm SL}(2,\Z)$, but for $n=3$ this is already quite a complicated space; for example the link of a rose is $2$-dimensional and homotopy equivalent to a wedge of eleven $2$-spheres.  

If $\G$ is a complete graph on $n$ vertices then $\AG$ is the free abelian group of rank $n$ and there are no $\G$-Whitehead partitions, so there are no $\G$-complexes other than the Salvetti (which is a torus), and  $\KG$ is a single point.    

If $\G$ has $3$ vertices $\{a,b,c\}$ and two edges, one joining $a$ to $c$ and one joining $b$ to $c$, then $\AG$ is the product  $F_2\times \Z=\langle a,b\rangle \times \langle c\rangle$ and the Salvetti is the product of a rose with a circle. Recall that a $\G$-Whitehead partition $\{P,P^*,{\rm lk}(P)\}$ is completely determined by giving either of its sides $P$ or $P^*$.
The only $\G$-Whitehead partitions in this case are those determined by  $P_1=\{a,b\}$ and $P_2=\{a,b^{-1}\}$, so the only $\G$ complexes are products of a circle with a rank $2$ graph (with no separating edges),  and $\KG$ can be identified with the spine $K_2$ of Outer space for the rank $2$ free group $\langle a,b\rangle$.  

If $\G$ has just a single edge from $a$ to $b$, then $\AG$ is the free product   $\Z^2at\Z=\langle a,b\rangle \ast \langle c \rangle$, and the Salvetti is a torus wedged with a circle. The only $\G$-Whitehead partitions are  those with sides   $P_1=\{a,c\}, P_2=\{b,c\}, P_3=\{a, c^{-1}\}$ and $P_4\{b,c^{-1}\}$.  The simplicial star in $\KG$ of a Salvetti is a square with corners $S^\PP,$ where $\PP=\{P_1,P_2\}, \{P_2,P_3\}, \{P_3,P_4\}$ or $\{P_4,P_1\}$. The  $\G$-complex associated to each corner is a torus cut into four squares, with an arc joining the two vertices.  There are  exactly four hyperplanes in this $\G$-complex which can be collapsed, so four stars of Salvettis fit together at each corner, and $\KG$ is a plane tiled by these stars.     The group  $\UAG$ is isomorphic to $\Z^2$ and acts on $\KG$ by translations.
 
 \subsection{Trees}
Recall that a RAAG $\AG$ is called {\em 2-dimensional} if $\G$ has no triangles, i.e. if the associated Salvetti complex is $2$-dimensional.  In \cite{BCV} we computed upper and lower bounds for the virtual cohomological dimension of $\Out$ for $2$-dimensional RAAGs.  For $\G$ a tree, these bounds agree, and we get 
$$\vcd(\Out)=e+2\ell-3$$
where $e$ is the number of edges and $\ell$ is the number of leaves.  

In the case that $\G$ is a tree, the only twists, $v \mapsto vw$, occur when $v$ is a leaf attached at $w$.  It follows that the twist subgroup $\TAG$ is free abelian of rank $\ell$.  It is easily checked that in this case, T$(\AG)$ is a normal subgroup of $\UAG$ and $\TAG \cap \UAG = \{1\}$, so $\Out$ is the semi-direct product of these two subgroups.  Since the cohomological dimension of $\TAG= \ell$, we conclude that the $\vcd (\UAG) \geq e +\ell - 3$.   In many (but not all) examples, the dimension of $\KG$ is exactly $e +\ell - 3$ and hence 
$\vcd (\UAG) = \dim (\KG)=  e +\ell - 3$. 

Consider, for example, the case in which $\G$ is a line made up of $e$ edges.  Label the vertices $v_0,v_1, \dots , v_e$.
Note that the generators of $\UAG$ are the partial conjugations by $v_i$,  $2 \leq i \leq e-2$,
and the folds $v_0 \mapsto v_0v_2^\pm$ and $v_e \mapsto v_ev_{e-2}^\pm$.  In particular, there are no $\G$-Whitehead 
partitions with $v_0, v_1, v_{e-1}$ or $v_e$ as maximal element.  

For $2 \leq i \leq e-2$, let $(R_i, v_i)$ be the $\G$-Whitehead pair corresponding to the partial conjugation by $v_i$ of all 
vertices to the right of $v_i$, that is $R_i=\{v_i,  v_{i+2}^{\pm 1}, \dots v_e^{\pm 1}\}$.  Similarly let $(L_i, v_i)$ be the pair corresponding to the partial conjugation by $v_i$ of all  vertices to the left of $v_i$.  Let $P_0=\{v_0, v_2\}$  and  $P_e = \{v_e, v_{e-2}\}$, so $(P_0, v_2)$ and $(P_e, v_{e-2})$ are the $\G$-Whitehead pairs corresponding to folds onto $v_0$ and $v_e$.  Then a typical maximal set of compatible 
partitions is of the form
$$ P_0, \, L_2 , \dots , L_j, \, R_{j+1}, \dots ,R_{e-2},\, P_e $$
Thus, the dimension of $\KG$ is $e-1=e + \ell -3$.

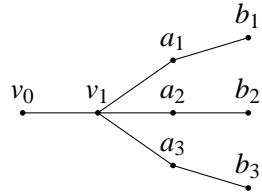
\begin{figure}
\begin{center}
\begin{tikzpicture}  
\draw [fill] (-1,0) circle [radius=0.03];
\draw [fill] (0,0) circle [radius=0.03];
\draw [fill] (1,0) circle [radius=0.03];
\draw [fill] (2,0) circle [radius=0.03];
\draw [fill] (1,.7) circle [radius=0.03];
\draw [fill] (2,1) circle [radius=0.03];
\draw [fill] (1,-.7) circle [radius=0.03];
\draw [fill] (2,-1) circle [radius=0.03];
\draw (-1,0) to  (2,0);
\draw (0,0) to   (1,.7) to (2,1);
\draw (0,0) to  (1,-.7) to (2,-1);
\node [above] at (-1,0) {$v_0$};
\node [above] at (0,0) {$v_1$};
\node [above] at (1,.7) {$a_1$};
\node [above] at (2,1) {$b_1$};
\node [above] at (1,0) {$a_2$};
\node [above] at (2,0) {$b_2$};
\node [above] at (1,-.7) {$a_3$};
\node [above] at (2,-1) {$b_3$};

 \end{tikzpicture}
\end{center}
\caption{An example for which $\dim(\KG) \neq e+\ell-3$}\label{counterexample}
\end{figure}

On the other hand, there exist examples of trees for which this equality does not hold.  We leave it as an exercise for the reader to show that for the tree $\G$ pictured in Figure~\ref{counterexample}, $\dim (\KG) =  e+\ell -2 = 9$.  Thus for this tree, we can only conclude that  $8 \leq \vcd(\UAG) \leq 9$.

\subsection{Strings of diamonds}   

\begin{figure}
\begin{center}
\begin{tikzpicture}
\draw (0,0)--(1,1)--(3,-1)--(5,1)--(6,0) -- (5,-1) -- (3,1) -- (1,-1)--(0,0);
\draw (6.5,.5)--(6,0)--(6.5,-.5);
\draw [fill] (6.5,0) circle [radius=0.05];
\draw [fill] (7,0) circle [radius=0.05];
\draw [fill] (7.5,0) circle [radius=0.05];
\draw (7.5,.5)--(8,0)--(7.5,-.5);
\draw (8,0)--(9,1)--(11,-1)--(12,0)--(11,1)--(9,-1)--(8,0);
\node [left] (c0) at (0,0) {$c_0$};
\draw [fill] (0,0) circle [radius=0.03];
\node [left] (c1) at (2,0) {$c_1$};
\draw [fill] (2,0) circle [radius=0.03];
\node [left] (c2) at (4,0) {$c_2$};
\draw [fill] (4,0) circle [radius=0.03];
\node [left] (c3) at (6,0) {$c_3$};
\draw [fill] (6,0) circle [radius=0.03];
\node [right] (cd1) at (8,0) {$c_{d-2}$};
\draw [fill] (8,0) circle [radius=0.03];
\node [right] (cd2) at (10,0) {$c_{d-1}$};
\draw [fill] (10,0) circle [radius=0.03];
\node [right] (e2) at (12,0) {$c_d$};
\draw [fill] (12,0) circle [radius=0.03];
\node [above] (t1) at (1,1) {$a_1$};
\draw [fill] (1,1) circle [radius=0.03];
\node [below] (b1) at (1,-1) {$b_1$};
\draw [fill] (1,-1) circle [radius=0.03];
\node [above] (t2) at (3,1) {$a_2$};
\draw [fill] (3,1) circle [radius=0.03];
\node [below] (b2) at (3,-1) {$b_2$};
\draw [fill] (3,-1) circle [radius=0.03];
\node [above] (t3) at (5,1) {$a_3$};
\draw [fill] (5,1) circle [radius=0.03];
\node [below] (b3) at (5,-1) {$b_3$};
\draw [fill] (5,-1) circle [radius=0.03];
\node [above] (td1) at (9,1) {$a_{d-1}$};
\draw [fill] (9,1) circle [radius=0.03];
\node [below] (bd1) at (9,-1) {$b_{d-1}$};
\draw [fill] (9,-1) circle [radius=0.03];
\node [above] (td) at (11,1) {$a_d$};
\draw [fill] (11,1) circle [radius=0.03];
\node [below] (bd) at (11,-1) {$b_d$};
\draw [fill] (11,-1) circle [radius=0.03];
\end{tikzpicture}
\end{center}
\caption{A string of $d$ diamonds}
\label{diamonds}
\end{figure}
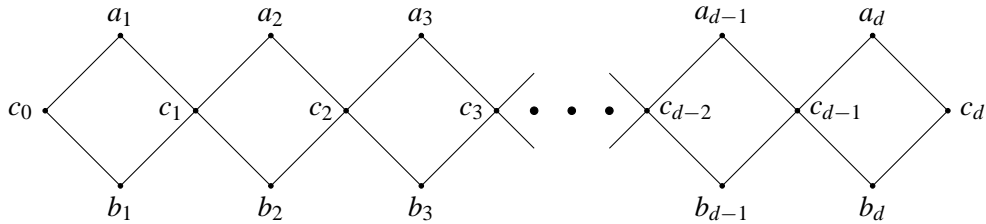

Let $\G$ be a ``string of $d$ diamonds," with $d\geq 3$ (see Figure~\ref{diamonds}).

This defines a $2$-dimensional $RAAG$ which has no twist automorphisms, so $\UAG=\Out$.  The upper and lower bounds from \cite{BCV} in this case are
$$3d+1\leq \vcd(\Out)\leq 6d-4.$$ 

We calculate the dimension of $\KG$ by finding the maximal possible size of a collection of compatible $\G$-Whitehead partitions.

Let $(P,m)$ be a $\G$-Whitehead pair with associated partition $\{P,P^\ast, {\rm lk}(m)\}$. 
\begin{itemize}
\item  There are no such pairs with $m=c_0$ or $m=c_d.$ 
\item Suppose $m=c_1.$  Removing the star of $c_1$ separates $\G$ into a left and right component, $L_1=\{c_0\}$ and $R_1$.
Since $c_0\leq c_1,$  either $c_0$ or $c_0^{-1}$ can be a singleton.  The possibilities for $P$ are $\{c_1,c_0\},$ $\{c_1^{-1},c_0^{-1}\},$   $\{c_1,c_0^{-1}\},$ $\{c_1^{-1},c_0\},$ $\{c_1,c_0,c_0^{-1}\}$ and  $\{c_1^{-1},c_0,c_0^{-1}\}.$  At most two of these are pairwise compatible, e.g. $\{c_1,c_0\}$ and  $\{c_1^{-1},c_0^{-1}\}.$  The case $m=c_{d-1}$ is symmetric.

\item Suppose $m=c_i$ with $2\leq i\leq d-2$.   Removing the star of $c_i$ again separates $\G$ into two components, $L_i$ and $R_i$.  Now there are no other possible singletons, so $P=\{c_i,L_i^\pm\}$ or $\{c_i,R_i^\pm\}$ (which are not compatible.)

\item If  $m=a_1$ the star separates $\G$ into $\{b_1\}$ and another component $C$; $b_1$ can be a singleton, and the possibilities for $P$ are
$\{a_ 1,b_1\}$ $\{a_ 1,b_1^{-1}\}$ and $\{a_1,b_1,b_1^{-1}\}$, $\{a_1, C^\pm\}$, $\{a_1, b_1, C^\pm\}$ and $\{a_1, b_1^{-1}, C^\pm\}$.  At most two of these are pairwise compatible, e.g. $\{a_ 1,b_1\}$ and $\{a_1, b_1, C^\pm\}= \{a_1^{-1}, b_1^{-1}\}^*$.  The case of $m=a_d$ is similar.

\item If  $m=a_i$   for $2\leq i\leq d-1$ the star separates $\G$ into three components, $K_i$ $\{b_i\}$ and $S_i$.  Possibilities for $P$ include $\{a_i,b_i\}, \{a_i,b_i, K_i^\pm\}$, and $\{a_i,b_i, K_i^\pm, S_i^\pm\}= \{a_i^{-1},b_i^{-1}\}^*,$ which is a maximal pairwise compatible collection with $m=a_i$.  The possibilities for $m= b_i$ are symmetric. 
  
\end{itemize}

One can check that the following sets  determine a collection of pairwise compatible $\G$-Whitehead partitions:

\begin{enumerate}
\item $ \{c_0,c_1\},\{c_{0}^{-1}, c_1^{-1}\}, \{c_{d-1},c_d\}, \{c_{d-1}^{-1}, c_d^{-1}\} $
\item  $\{c_i,L_i^\pm\} \hbox{ for } 2\leq i\leq d-2$
\item $\{a_i,b_i\}, \{a_i^{-1},b_i^{-1}\}, \hbox{ for } 1\leq i\leq d.$
\item $\{a_i,b_i,K_i\}\hbox{ for } 2\leq i\leq d-1$ 
\end{enumerate}

By the remarks above, this is a maximal collection.  There are $4+(d-3)+2d+(d-2)=4d-1$ of them, so $\dim(\KG)=4d-1$.

\subsection{Butterflies} 

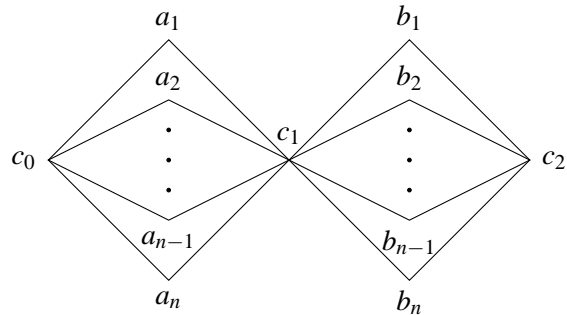
\begin{figure}
\begin{center}
\begin{tikzpicture}[scale=.8]
\draw (0,0)--(2,2)--(4,0)--(2,-2)--(-2,2)--(-4,0)--(-2,-2)--(0,0);
\draw (0,0)--(2,1)--(4,0)--(2,-1)--(-2,1)--(-4,0)--(-2,-1)--(0,0);
 \draw [fill] (2,.5) circle [radius=0.03];
  \draw [fill] (2,0) circle [radius=0.03];
   \draw [fill] (2,-.5) circle [radius=0.03];
 \draw [fill] (-2,.5) circle [radius=0.03];
  \draw [fill] (-2,0) circle [radius=0.03];
   \draw [fill] (-2,-.5) circle [radius=0.03];
   \node [above] (a1) at (-2,2) {$a_1$};
      \node [above] (a2) at (-2,1) {$a_2$};
         \node [below] (an1) at (-2,-1) {$a_{n-1}$};
         \node [below] (an) at (-2,-2) {$a_n$};
           \node [above] (b1) at (2,2) {$b_1$};
      \node [above] (b2) at (2,1) {$b_2$};
         \node [below] (bn1) at (2,-1) {$b_{n-1}$};
         \node [below] (bn) at (2,-2) {$b_n$};
   \node[left] (l) at (-4,0) {$c_0$};
   \node[right] (r) at (4,0) {$c_2$};
   \node [above] (c) at (0,.1) {$c_1$};
\end{tikzpicture}
\end{center}
\caption{Butterfly graph}
\label{butterfly}
\end{figure}

Let $\G$ be the ``butterfly" graph shown in Figure~\ref{butterfly} (which is a string of two diamonds if $n=2$).
 Again $\AG$ is $2$-dimensional with no twists, so $\UAG=\Out$.  The bounds on the $\vcd$ from \cite{BCV} in this case are 
 $$4n-1\leq \vcd(\Out)\leq 4n.$$
 The maximal size of a compatible set of $\G$-Whitehead partitions is $4n-1$, realized for example by the partitions determined by the following sets $P$:
\begin{enumerate}
\item $ \{c_1,c_0\},\{c_1^{-1}, c_0^{-1}\}, \{c_1, c_0, c_2\} $
\item $\{a_1,a_2,\ldots a_i\}, \{a_1^{-1},a_2^{-1},\ldots a_i^{-1}\}, \hbox{ for } 2\leq i\leq n.$
\item $\{b_1,b_2,\ldots b_i\}, \{b_1^{-1},b_2^{-1},\ldots b_i^{-1}\}, \hbox{ for } 2\leq i\leq n.$
\end{enumerate}
Thus the dimension of $\KG$ is $4n-1$. Since this matches the lower bound this is equal to the $\vcd$ of $\Out$.

\section{Contractibility of the  \texorpdfstring{$\boldsymbol\G$}{G}-spine}

Our strategy for showing that   $\KG$ is contractible is to view it as the union of stars of \roses.  We first define a norm which totally orders the \roses.  We then construct the space by starting with the star of the \rose\   of minimal norm, then attaching the stars of the rest of the \roses\  in increasing order.  We check at each stage that we are attaching along a contractible subcomplex.

\subsection{The norm of a \rose.}

The norm is defined using lengths of conjugacy classes of elements of $\AG$, and we begin with some observations about these lengths.  Let $\sigma=(\SG,\alpha)$ be a \rose.  For any conjugacy class  $g$, define $\ell_\sigma(g)$ to be the minimal length of a  word $\word$ in the free group $F(V)$ representing an element of the conjugacy class of $\alpha^{-1}(g)$ in $\AG$.    In particular, if $\sigma=(\SG,id)$, then $\ell_\sigma(g)$ is the minimal word length of an element of $g$.  

Normal form for elements of $\AG$ (see, e.g., \cite{Cha}) implies that $\ell_\sigma(g)$ is well-defined.  Since the vertices $V$ of $\G$ can be identified with  the edges in the 1-skeleton of $\SG$,   $\ell_\sigma(g)$ can also be thought of as the length of a minimal edge-path in the 1-skeleton of $\SG$ representing $\alpha^{-1}(g)$.  If $\alpha$ is an isometry of $\SG$, then this is the same as the length of $g$, reflecting the fact that $\sigma=(\SG,\alpha)$ is equal to $(\SG,id)$ as  a \rose.

Let $\Gcal=(g_1,g_2,\ldots)$ be a list of all conjugacy classes in  $\AG$, and let $\Gcal_0$ be the set of conjugacy classes which can be represented by words of length at most $2$.  

\begin{definition}For  a \rose\  $\sigma=(X,\alpha),$  the {\em norm} $\|\sigma\|=(\|\sigma\|_0,\|\sigma\|_1,\|\sigma\|_2, \ldots) \in \Z\times\Z^\Gcal$ is defined by
$$
     \|\sigma\|_0=\sum_{g\in\Gcal_0} \ell_\sigma(g), \qquad
     \|\sigma\|_i=\ell_\sigma(g_i) \hbox{ for } i\geq 1.
$$    
We consider $\Z\times\Z^\Gcal$ as an ordered abelian group, with the lexicographical ordering.  Denote the identity element by $\0 = (0,0, \dots )$.  We say an element $x=(x_0,x_1,\ldots)\in \Z\times\Z^\Gcal$ is \emph{negative} if $x<\0$, and \emph{strongly negative} if its first coordinate $x_0$ is negative.
\end{definition}

 \begin{lemma}\label{unique} The \rose\  $(\SG,id)$ is the unique \rose\ of minimal norm.  Indeed, for any other \rose\  $\sigma$,   $\|(\SG,id)\|_0 < \|\sigma\|_0$.
 \end{lemma} 

\begin{proof} An easy calculation shows that if $|V|=m$, the first coordinate   $\|(\SG,id)\|_0$ is equal to  
$2m^2+2m,$  and that is the minimal possible value. It suffices to show that the first coordinate of the norm of any other  \rose\ is strictly larger.

Suppose $\sigma=(\SG,\alpha)$ is another marked Salvetti complex with $\|\sigma\|_0 =2m^2+2m$. 
Then $\alpha$ must permute the conjugacy classes in $\Gcal_0$ since otherwise, $\ell_\sigma(g) > 2$ for some $g \in \Gcal_0$.  In fact, a stronger statement holds:  $\alpha$ must permute the conjugacy classes of $V^\pm\subset \Gcal_0$, since if $\ell_\sigma (v) =2$ for some $v \in V$, then 
$\ell_\sigma (v^2) = 4$.  
Thus, $\alpha$ induces a permutation of the directed edges of $\SG$.  Moreover, if two edges of $\SG$ span a cube, then their images must also span a cube, since if $v,w \in V$ do not commute, then neither do any conjugates of $v$ and $w$.  Thus after composing with an isometry of $\SG$ we may assume $\alpha$ takes every element of $V$ to a conjugate of itself.   

Let $V=\{v_1, \dots v_m\}$ and choose an automorphism $\alpha_1 \in Aut(\AG)$ representing $\alpha$ such that $\alpha_1(v_1)=v_1$. 
Say $\alpha_1^{-1}(v_2)= av_2a^{-1}$ where $a$ is of minimal length (i.e., $av_2a^{-1}$ is a reduced word).  Then $\ell_\sigma(v_1v_2) =2$ implies that the cyclic reduction of $v_1av_2a^{-1}$ is a word of length 2.  The only way this can happen is if $a$ lies in the centralizer of $v_1$.  Thus, we can compose $\alpha_1$ with conjugation by $a^{-1}$ to get a new representative $\alpha_2$ which acts as the identity on both $v_1$ and $v_2$. 

Now repeat with $v_3$.  Say $\alpha_2^{-1}(v_3)=bv_3b^{-1}$.  Arguing as above, $b$ must lie in the intersection of the centralizers $C(v_1) \cap C(v_2)$, so composing $\alpha_2$ with conjugation by $b^{-1}$ gives a representative for $\alpha$ which acts as the identity on $v_1,v_2$ and $v_3$.  Continuing in this manner, we see that $\alpha$ has a representative which is the identity on all of $V$, that is, $\alpha$ is homotopic to the identity map. 

\end{proof}

\begin{corollary}\label{claimtwo}  
 Given a \rose\  $\sigma=(\SG,\alpha)$ there is a finite set of conjugacy classes  $\Gcal_\sigma \subset \AG$ such that 
 $\sigma$ is uniquely determined by $\sum_{g\in\Gcal_\sigma}\ell_\sigma(g)$.   
 \end{corollary}

 \begin{proof}   
Replace  $\Gcal_0$ in the proof of Lemma~\ref{unique} by the set of $g$ with $\ell_\sigma(g) \leq 2$.  
 \end{proof}

In particular, no two \roses\ have the same norm.  In section~\ref{sec:reductive} we will show that the norm induces a well-ordering of the set of  \roses\ but we need some preparation first.

\subsection{Effect of a Whitehead move on the norm}\label{moves}

We want to construct $\KG$ by adding stars of \roses\ to the star of $(\SG, id)$ in order of increasing norm.  Since adjacent stars are connected by \GW moves, we will need to understand how lengths of conjugacy classes change under these moves. We will do this first by using the geometric interpretation of $\ell_\sigma(g)$.  

Let $\sigma=(\SG,\alpha)$ be a  \rose, $\Pbold$ a \GW partition of $V^\pm$ and $v\in {\rm max}(P)$.
The length  $\ell_\sigma(g)$ is the length of the shortest edge path in the 1-skeleton of $\SG$ representing the free homotopy class of $\alpha^{-1}(g)$.  To understand what happens to this length under the Whitehead automorphism  $(\Pbold, v)$ we will find a minimal length edge path in $\SGP$ representing $\alpha^{-1}(g)$ and apply Lemma \ref{minimal paths}.

 Recall that $\SGP$ has exactly two vertices, corresponding to the two sides of $\Pbold$, and an edge labelled $e_\Pbold$ between them.  For $u \in {\rm lk}(\Pbold)$,  there is a loop labelled $e_u$ at each vertex.  For $u \notin {\rm lk}(\Pbold)$, there is a unique edge labelled $e_u$ with initial vertex corresponding to the side of $\Pbold$ containing $u$ and terminal vertex corresponding to the side containing $u^{-1}$.  Let $\word$ be a  cyclically reduced word for $\alpha^{-1}(g)$.  Identifying directed edges of $\SG$ with $V^\pm$,  let $p$ be the edge path  in $\SG$ labelled by $\word$.  We can lift $p$ to a loop $\tilde p$ in $\SGP$ as follows. 
  If the support of $\word$ lies entirely in ${\rm lk}(\Pbold)$, $p$ lifts to an edge path $\tilde p$ of the same length at either vertex.  Otherwise, cyclically permuting $\word$ if necessary, let $\word=u_1 \dots u_k$ where $u_i \in V^\pm$ and $u_1 \notin {\rm lk}(\Pbold)$.  Then $u_1$ corresponds to a unique directed edge $e_1$ in $\SGP$.  If $u_1$ and $u_2^{-1}$ both lie in $\overline P=P \cup {\rm lk}(P)$ or both in $\overline P^*= P^*  \cup {\rm lk}(P)$, then 
  $u_2$ lifts to a directed edge $e_2$ whose initial vertex equals the terminal vertex of $e_1$.  Hence $u_1u_2$ lifts to the path $e_1e_2$.  If not, insert the edge $e_\Pbold$ (appropriately oriented) to get a path $e_1e_\Pbold e_2$ which projects to $u_1u_2$.  Now repeat this process with each $u_i$ to obtain the loop $\tilde p$. 

It is easy to see that $\tilde p$ is a minimal length lift of $p$.  To see that it is a minimal length representative for  
$\alpha^{-1}(g)$, note that any other minimal word $\word' $ for $\alpha^{-1}(g)$ can be obtained from $\word$ by interchanging commuting pairs $u_iu_{i+1}$.  But for such a pair, the edges $e_{u_i}$ and $e_{u_{i+1}}$ span a square in $\SGP$, so they can be traversed in either order without crossing $e_\Pbold$. 
It follows that the length of $\tilde p$ is independent of choice of $\word$.   

To keep track of the lengths of these paths, we introduce some new notation. Set
\begin{itemize}
\item $|\Pbold |_{\word}=$ the number of times $\tilde p$ traverses the edge $e_\Pbold$, or equivalently, the number of (cyclically) adjacent letters $u_iu_{i+1}$ in $\word$ such that $u_i$ and $u_{i+1}^{-1}$ do not both lie in  $\overline P$ or both in $\overline P^*$,
\item  $|v|_{\word}=$ the number of occurrences of  $v$ or $v^{-1}$ in $\word$. 
\end{itemize}

\begin{lemma} \label{norm1}  Let $\sigma=(\SG,\alpha)$ be a \rose, 
let $\varphi=(\Pbold,v)$ be a Whitehead automorphism, let $g$ be a conjugacy class in $\AG$ and let $\word$ be a minimal length word representing $\alpha^{-1}(g)$.  Then
$$\ell_{\sigma^\Pbold_v}(g) =  \ell_\sigma(g) +  |\Pbold|_{\word}   - |v|_{\word}.$$
More generally, if $\sigma'$ is obtained from $\sigma$ by blowing up a compatible collection of  \GW partitions 
$\PP=\{\Pbold_1,\ldots,\Pbold_k\}$ and collapsing  a tree-like set of hyperplanes $\Hcal=\{H_1, \dots ,H_k\}$ dual to edges labelled $e_{v_i}$,  then
$$\ell_{\sigma'}(g) =  \ell_\sigma(g) + \sum_{i=1}^k  |\Pbold_i|_{\word}   - \sum_{i=1}^k  |v_i|_{\word} .$$
\end{lemma} 

\begin{proof}  First consider the blow-up-collapse for a single Whitehead pair $(\Pbold,v)$.    Let $p$ and $\tilde p$ be as above.  By construction,  $\ell_\sigma(g)= length(p) = length (\tilde p) -  |\Pbold|_{\word}$.
Collapsing the hyperplane in $\SGP$ dual to the edge labelled $e_v$ gives the marked Savetti $\sigma^\Pbold_v=(\SG, \varphi \alpha)$.    Let $p'$ be the image of $\tilde p$ under this collapse.  By Lemma \ref{minimal paths}, $p'$ is a minimal length representative for $(\varphi \alpha)^{-1}(g)$.  Hence
$\ell_{\sigma^\Pbold_v}(g) =length(p') = length (\tilde p) - |v|_{\word}$.  This proves the first statement.

For the second statement,  let $c \colon \SGPP \to \SG$ be the canonical projection and let $c_{\Hcal}$ be the collapsing map onto $\SG^{\PP}_{\Hcal}$.  Let $\tilde p$ be a minimal edge path in $\SGPP$ representing $\alpha^{-1}(g)$.  (Here we identify the fundamental group of $\SGPP$ and $\SG$ via $c$.)
Let $p=c(\tilde p)$ and  $p'=c_\Hcal(\tilde p)$.
Then by Lemma \ref{minimal paths},  $p$ and $p'$ are minimal paths in their homotopy class.  In particular, $p$ corresponds to a minimal word $\word$ representing $\alpha^{-1}(g)$ so the number of edges of $p$ (and hence also of $\tilde p$) labelled $e_{v_i}$ equals $|v_i|_{\word}$.  

Collapsing $\SGPP$ to a single blow-up $\SG^{\Pbold_i}$ maps $\tilde p$ to a minimal lift $p_i$ of $p$, hence by the discussion above, the number of edges of $p_i$ (and hence also of $\tilde p$) labelled $e_{\Pbold_i}$ equals $|\Pbold_i|_{\word}$.  It now follows that 
$$\ell_{\sigma}(g) =  length(p) = length(\tilde p) - \sum |\Pbold_i|_{\word} $$
$$\ell_{\sigma'}(g) =  length(p') = length(\tilde p) - \sum |v_i|_{\word}$$
\end{proof}

\begin{remark}  The hypothesis that every hyperplane in $\Hcal$ be dual to an edge labelled $e_{v_i}$ is crucial in this lemma.  In general, a tree-like set $\Hcal=\{H_1, \dots , H_k\}$ in $\SGPP$ may contain hyperplanes dual to edges labelled $e_{\Pbold_i}$.  Collapsing these hyperplanes first to get a smaller $\G$-complex, we see that $\sigma'=
\sigma^{\PP}_{\Hcal}$ can be obtained from $\sigma$ by a blow-up-collapse satisfying the conditions of the lemma. 
\end{remark}

We put all of this information together for a \rose\   $\sigma=(\SG,\alpha)$ and a Whitehead automorphism $\varphi=(P,v)$ by defining absolute values $|P|_\sigma$ and $|v|_\sigma $ in $\Z\times\Z^\Gcal$ coordinate-wise, i.e.   
$$|P|_\sigma= (|P|_{0},|P|_{\word_1}, |P|_{\word_2},\ldots),$$  where
\begin{itemize}
 \item $\word_i$ is a minimal length word representing the conjugacy classe $\alpha^{-1}(g_i),$   and 
 \item $|P|_{0}=\sum_{\word\in \mathcal W_0}|P|_\word$ for a set of words $\mathcal W_0$ representing the $\alpha^{-1}(g)$ for $g\in\Gcal_0$.
 \end{itemize}
 Similarly, define
$$|v|_\sigma=(|v|_{0},|v|_{\word_1},|v|_{\word_2},\ldots),$$ 
where $|v|_{0} =\sum_{\word\in \mathcal W_0}|v|_\word$.
 
Lemma~\ref{norm1} can now be restated as
     
 \begin{corollary}\label{norm2} Let $\PP, \Hcal$ be as in Lemma \ref{norm1}.  
Then $$\|\sigma^\PP_\Hcal\| = \|\sigma\| + \sum  |\Pbold_i|_{\sigma}   - \sum |v_i|_{\sigma}. $$
\end{corollary}

\begin{definition} A \GW partition $\Pbold$ is   \emph{reductive} for a \rose\  $\sigma$   if for some $v\in {\rm max}(P)$ the Whitehead automorphism $\varphi=(P,v)$  reduces $\|\sigma\|,$ that is,  $\|\sigma^\Pbold_v\|< \|\sigma\|$, or equivalently, $|\Pbold_\sigma|  < |v_{\sigma}|$.  It is \emph{strongly reductive} if the first coordinate   $\|\sigma^\Pbold_v\|_0$  is less than $\|\sigma\|_0.$  \end{definition}

By Corollary~\ref{claimtwo}, $\sigma$ and $\sigma^\Pbold_v$ cannot have the same norm, since they are different \roses.

\begin{corollary}\label{reductive} Let $\PP, \Hcal$ be as in Lemma \ref{norm1}.   If $\|\sigma^\PP_\Hcal\| < \|\sigma\|$, then some $\Pbold_i\in\PP $ is reductive for $\sigma$.  If $\sigma^\PP_\Hcal$ is strongly reductive, then so is some $\Pbold_i$.
\end{corollary}

\begin{proof}  By Theorem \ref{swap},  the elements of $\Hcal$ can ordered such that if $e_{v_i}$ is the edge dual to $H_i$, then  $(P_i,v_i)$ is a \GW  pair.  If $\|\sigma^\PP_\Hcal\| < \|\sigma\|$, then by the previous corollary, 
 $|\Pbold_i|_{\sigma}   -  |v_i|_{\sigma} < \0$ for some $i$, so 
 $ \|\sigma^{\Pbold_i}_{v_i}\| < \|\sigma\| $.  The same argument applied to the first coordinate of the norm shows that if $\sigma^\PP_\Hcal$ is strongly reductive, then so is some $\sigma^{\Pbold_i}_{v_i}$.

\end{proof}

\subsection{Star graphs}

There is a convenient combinatorial way to keep track of computations such as those we did in  Section~\ref{moves} using a diagram called a {\it star graph}. 
Star graphs have been extensively used  to study free groups and their automorphisms (see, e.g., \cite{Hoare}).  The star graph $\St(\word)$ of a cyclically reduced word $\word$ in a free group $F(V)$ is defined by taking a vertex for each element of $V^{\pm}$ and  an edge from $x$ to $y$ for every occurrence of $xy^{-1}$ as a (cyclic) subword of $\word$.   If we consider $F(V)$ as a right-angled Artin group on the discrete graph $\G$, then  $|v|_\word$ (as defined in the previous section) is equal to the valence  of a vertex $v$ in $\St(\word)$,  and for
any partition $\Pbold$ of $V^{\pm}$, $|\Pbold|_\word$ is equal to
the number of edges in the star graph with one vertex in  $P$ and one vertex in $P^*$.  
Since the star graph $\St(\word)$ depends only on $\word$ in the case of a free group, it can be used to compute $|\Pbold|_\word$ for any $\Pbold$.

We would like to imitate this construction for more general $\AG$, but for a conjugacy class $g$ of $\AG$ and \GW partition $\Pbold$, to compute $|\Pbold|_{\word}$,  we need  to count how many times a minimal path $\tilde p$ in the $\G$-complex crosses an edge labelled $e_\Pbold$.  This involves counting not only when $\word$  crosses from $P$ to $P^*$, but when it is forced to cross from $\overline P= P \cup {\rm lk}(P)$ to $\overline P^*=P^* \cup {\rm lk}(P^*)$.  Since our star graphs must  take into account the link of $\Pbold$ they cannot be defined independently of the partition. 

Consequently, for a \GW partition $\Pbold=\{P,{\rm lk}(P),P^*\}$ and a cyclically reduced word $\word=u_1 \dots u_k$ we define the star graph $\St_\Pbold(\word)$ as follows.  The vertices of $\St_\Pbold(\word)$ are the elements of the disjoint union of  $\overline P=P\cup {\rm lk}(P)$ and $\overline P^*=P\cup {\rm lk}(P)$, i.e., we have two copies of ${\rm lk}(P)$ instead of one.   View $\word$ as a cyclic word and set $u_{k+1}=u_1$.  Beginning with $i=1$  draw an edge from  $u_i$ to $u_{i+1}^{-1}$ staying within $\overline P$ or $\overline P^*$ whenever possible. If every $u_i$ lies in ${\rm lk}(P)$, then the star graph can be drawn entirely in $\overline P$ (or in $\overline P^*$).  Otherwise, we may cyclically permute $\word$ so that $u_1$ does not lie in ${\rm lk}(P)$, in which case there is no choice of where to start.  See figure \ref{stargraph} for an example.

The quantities $|v|_\word$,  for $v \notin {\rm lk}(P)$,  and $|\Pbold|_{\word}$ can now be read off the star graph  $\St_\Pbold(\word)$.  Namely,  $|v|_{\word}$ is equal to the valence of the vertex $v$, while $|\Pbold|_{\word}$ equals the number of edges with one vertex in $\overline P$ and one vertex in $\overline P^*$.  

We will need to compare star graphs for the same word with respect to different partitions, but  the graph we have constructed depends on the partition $\Pbold$, not just on the word $\word$. 
To solve this problem, we will need to consider slightly more general decompositions of $V^{\pm}$, and a more general definition of a star graph.

Fix a symmetric subset $L \subset V^\pm$, a decomposition  $A_1+\ldots+A_k $ of the complement $L^c$.
The  {\em  star graph} 
$\St^L_{A_1,\ldots,A_k}(\word)$  is constructed as follows.  Take a copy $L_i$ of $L$ for each $A_i$ and let $\overline A_i = A_i \cup L_i$.  The vertices of  $\St^L_{A_1,\ldots,A_k}(\word)$ are the elements of the (disjoint) union  of the  $\overline A_i$.   We draw  the star graph $\St^L_{A_1,\ldots,A_k}(\word)$ by first drawing circles to isolate each $\overline A_i$.  The idea is then to draw the edges of $\St^L_{A_1,\ldots,A_k}(\word)$ in order, avoiding crossing circles whenever possible. 

More precisely, we proceed as follows.  If all letters of $\word$ are in $L$, we will draw the entire star graph with vertices in $L_1.$   Otherwise, list all of the 2-letter subwords $xy$ of $\word$ in order (cyclically), starting at a letter $x\in L^c$.  Since  $\{A_i\}$ partitions $L^c$, $x$ lies in a unique $A_i$.  If $y\in L^c$, there are unique vertices labelled $y$ and $y^{-1}$, so we have no choice: we draw an edge from $x$ to $y^{-1}$ and start the next edge at $y$.  If $y\in L$, draw an edge from $x$ to the copy of $y^{-1}\in L_i$, and start the next edge at $y\in L_i$.  We continue in this way, remaining inside each $\overline A_j$-circle as long as possible.   
Note that if $\Pbold$ is a \GW partition and $L = {\rm lk}(\Pbold)$, then  $\St^L_{P,P^*}(\word)$ is precisely the graph $\St_\Pbold(\word)$ constructed above.

\begin{figure}
\begin{center}
 \begin{tikzpicture}[scale=.4]
\draw [rounded corners] (-1,-1) rectangle (6,12);
\draw [rounded corners] (0,0) rectangle (5,4);
\node (lk) at (2.5,1) {${\rm lk}(P)$};
\node (y1label) at (1,3) {$y$};
\coordinate (y1) at (1.5,2.5);
\fill  (y1) circle(.1cm);
\node (ybar10) at (4,3) {$y^{-1}$};
\coordinate  (ybar1) at (3,2.5);
\fill  (ybar1) circle(.1cm);
\draw [rounded corners] (0,5) rectangle (5,11);
\node (P) at (1,10) {$P$};
\node (x0) at (4,7) {$x$};
\coordinate (x) at (3.5,6.5);
\fill (x) circle(.1cm);
\node (u) at (1,7) {$u$};
\fill  (1.5,6.5) circle(.1cm);
\node (v) at (3,9) {$v$};
\fill  (2.5,8.5) circle(.1cm);
\begin{scope}[xshift=8cm]
\draw [rounded corners] (-1,-1) rectangle (6,12);
\draw [rounded corners] (0,0) rectangle (5,4);
\node (lk) at (2.5,1) {${\rm lk}(P)$};
\node (y20) at (1,3) {$y$};
\coordinate (y2) at (1.5,2.5);
\fill  (y2) circle(.1cm);
\node (ybar20) at (4,3) {$y^{-1}$};
\coordinate (ybar2) at  (3,2.5);
\fill  (ybar2) circle(.1cm);
\draw [rounded corners] (0,5) rectangle (5,11);
\node (Pstar) at (1,10) {$P^*$};
\node (z0) at (4,7) {$z$};
\coordinate (z) at (3.5,6.5);
\fill  (z) circle(.1cm);
\node (xinv) at (1,7) {$x^{-1}$};
\coordinate (xinv) at (1.5,6.5);
\fill  (xinv) circle(.1cm);
\node (zinv0) at (3.25,9.25) {$z^{-1}$};
\coordinate (zinv) at (2.5,8.5);
\fill  (zinv) circle(.1cm);
\end{scope}
\draw [red, thick] (y1) to (x) to (zinv);
\draw[red, thick] (ybar1) to (xinv) to (y2);
\draw [red, thick] (ybar2) to (z);
 \end{tikzpicture}
\end{center}
\caption{$\St_\Pbold(xy^{-1}xzy)=\St^{{\rm lk}(P)}_{P,P*}(xy^{-1}xzy)$}
\label{stargraph}
\end{figure}  

If $v\in A_i$,   the valence of $v$ in $\St^L_{A_1,\ldots,A_k}(\word)$  is equal to the number of occurrences of $v$ or $v^{-1}$ in $\word$, and if $v\in L$, then the number of such occurrences is equal to the sum of the valences of the copies of $v$ in the $L_i$.

\subsection{Counting lemmas} This section contains several elementary counting lemmas related to star graphs which are at the heart of the proofs in the next section. 

As above, let $L$ be a symmetric subset of $V^\pm$.  For a subset $A \subset L^c,$ denote 
$$A^*= A^c \setminus L, \qquad \overline A= A \cup L_A,  \qquad  \overline A^*= A^* \cup L_{A^*}$$
 where  $L_A$ and $L_{A^*}$ are copies of $L$.

\begin{definition} For a  cyclically reduced word in $\word\in F(V)$  and  disjoint subsets $A$ and $B$ of $L^c$,  define  the {\em dot product}  $(A.B)^L_{\word}$ to be the number of edges of $\St^L_{A,B,(A+B)^*}(\word)$ with one vertex in $\overline A$ and one vertex in $\overline B$.
\end{definition}

The dot product $(A.B)^L_{\word}$ can also be described  as the number of cyclic subwords of  $\word$ of the form $a{\bf u}b^{-1}$ or $b{\bf u}a^{-1}$ for $a\in A, b\in B$ and ${\bf u}$  a word in $L$.  If $B=A^*$, then $A+B=L^c$,
so $(A+B)^* = \emptyset.$   In this case, no edge of the star graph enters the $\overline{(A+B)}^*$-circle, so for the purposes of our computations, we can identify  
$\St^L_{A,B,(A+B)^*}(\word)$ with $\St^L_{A,A^*}(\word).$

\begin{definition}
For a cyclically reduced word $\word$  and a subset $A\subset L^c$,  define the {\em absolute value} of $A$ by $|A|^L_{\word}=(A. A^*)^L_{\word}=$  the number of edges of $\St^L_{A,A^*}(\word)$ with one vertex in $\overline A$ and one vertex in $\overline A^*$.   
\end{definition}

\begin{example}  Let $\Pbold=\{P,P^*,{\rm lk}(P)\}$ be a \GW partition.  Then  $(P.P^*)^{{\rm lk}(P)}_{\word}=|P|^{{\rm lk}(P)}_{\word}=|P^*|^{{\rm lk}(P)}_{\word}=|\Pbold|_\word$.   
\end{example}

Our justification for calling $(A.B)^L_{\word}$ a ``dot product" rests partly on the observation that  $(A.B)^L_{\word}=(B.A)^L_{\word}$. We also have the following linearity relation.

\begin{lemma}\label{disjoint} Let $L\subset V^\pm$ be a symmetric subset, and let 
$A,B$ and $C$ be disjoint subsets of $L^c$.  Then $(A.(B+C))^L_{\word}=(A.B)^L_{\word}+(A.C)^L_{\word}$.
\end{lemma}

\begin{proof}  The  star graph $\St^L_{A,B+C, (A+B+C)^*}(\word)$ is obtained from $\St^L_{A,B,C,(A+B+C)^*}(\word)$ by identifying the vertices of $L_B$ and $L_C$.  
 Thus the number of edges between $\overline A$ and   $\overline{(B+C)}$ is equal to the number between $\overline A$ and $\overline B$ plus the number between $\overline A$ and $\overline C$, as desired.    
\end{proof}

In the following three lemmas we fix $\word$ and omit it from the notation for simplicity.   

\begin{lemma}\label{AB} Let $L$ be a symmetric subset of $V^\pm$.  For any subsets $A$ and $B$ of  $L^c$, 
$$|A\cap B^*|^L+|A^*\cap B|^L= |A|^L+|B|^L-2((A\cap B).(A^*\cap B^*))^L.$$ 
\end{lemma}
\begin{proof} We set notation according to the figure below:
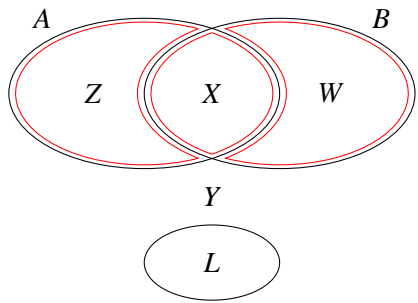
\begin{center}
 \begin{tikzpicture} [xscale=.45, yscale=.25]
  \draw (0,0) circle (4cm); 
  \node (A) at (-3,4) {$A$};
   \draw [red] (65:3.8) arc (65:295:3.8cm);
    \draw [red] (58.5:3.8) arc (58.5:-58.5:3.8cm);
    \draw [red] (55.5:4.2) arc (55.5:-55.5:4.2cm);
    
  \begin{scope}[xshift=4cm]
   \draw (0,0) circle (4cm);
   \node (B) at (3,4) {$B$}; 
   \draw [red] (122:3.8) arc (122:238:3.8cm);
    \draw [red] (115:3.8) arc (115:-115:3.8cm);
      \draw [red] (125:4.2) arc (125:235:4.2cm);
    \end{scope}
    
    \node (X) at (2,0) {$X$};
     \node (Y) at (2,-5.5) {$Y$};
      \node (Z) at (-1.5,0) {$Z$};
       \node (W) at (5.5,0) {$W$};
    \draw (2,-9) circle (2cm);
    \node (L) at (2,-9) {$L$};
 \end{tikzpicture}
\end{center}
Thus $X=A\cap B$, $Z=A\cap B^*$, $W= A^*\cap B$ and $Y= A^* \cap B^*$ and we are trying to show $$|Z|^L+|W|^L=|A|^L+|B|^L-2(X.Y)^L$$
 We calculate
$$
\begin{array}{rcl}
 |A|^L& = &((X+Z).(Y+W))^L  = (X.Y)^L+(X.W)^L+ (Z.Y)^L + (Z.W)^L\\
|B|^L & = &((X+W).(Y+Z))^L =  (X.Y)^L + (Y.W)^L + (Z.X)^L + (Z.W)^L\\
|Z|^L & = &(Z.(X+Y+W))^L = (Z.X)^L + (Z.Y)^L + (Z.W)^L\\
 |W|^L & =& (W.(X+Z+Y)^L = (W.X)^L+(W.Z)^L+(W.Y)^L.
\end{array}
$$
 The result follows.
\end{proof}

\begin{lemma}\label{difference} Let $L_0\subset L$ be symmetric subsets of $V^\pm$  and let $A\subset C \subset L^c$.  Then $ |A|^{L_0} - |A|^{L}  \leq  |C|^{L_0} - |C|^{L}.$
\end{lemma}

\begin{proof} $|A|^{L}$ counts the number of subwords of the form $a.{\bf u}.b^{-1}$ or $b.{\bf u}.a^{-1},$ with $a\in A$, $b\in (A+L)^c$ and ${\bf u}$ a (possibly empty) word in elements of $L$. Notice that each such subword also contributes exactly one to $|A|^{L_0}$. 

Let $\mathcal D$ be the set of all words in elements of $L$ that use at least one letter which is not in $L_0$.  The only other contributions to  $|A|^{L_0}$  come from subwords of the form $a^\prime.{\bf u}.a^{-1}$ for $a,a^\prime\in A$ and ${\bf u}\in \mathcal D$; each of these subwords contributes 2.  Thus the difference $ |A|^{L_0} - |A|^{L}$ is the number of subwords $a^\prime.{\bf u}.a^{-1}$ with ${\bf u}\in\mathcal D$.

We now do the same computation for $C$.  Since $A\subset C$, there are at least as many words of the form $c^\prime.{\bf u}.c^{-1}$ with ${\bf u}\in \mathcal D$  as words $a^\prime.{\bf u}.a^{-1}$, and the lemma is proved.
\end{proof}

\begin{lemma}\label{ABL}  Let $L_1, L_2$ be a symmetric subsets of $V^\pm$ and let  $L=L_1 \cup L_2.$
Then for any subsets $A$ and $B$ of  $L^c,$ 
$$|A\cap B^*|^{L_1}+|A^*\cap B|^{L_2}  \leq  |A|^{L_1}+|B|^{L_2}.$$ 
\end{lemma}

\begin{proof}   By Lemma \ref{difference} applied to $(A\cap B^*) \subset A \subset L^c$, we have
$$|A\cap B^*|^{L_1} + |A|^L \leq   |A\cap B^*|^{L} +|A|^{L_1}. $$
By Lemma \ref{difference} applied to $(A^*\cap B) \subset B \subset L^c$,
$$|A^*\cap B|^{L_2} + |B|^L \leq  |A^*\cap B|^{L} + |B|^{L_2}. $$
Adding these and applying Lemma \ref{AB} gives
\begin{align*}
|A\cap B^*|^{L_1}+|A^*\cap B|^{L_2} + |A|^L + |B|^L
  &\leq |A\cap B^*|^{L}+|A^*\cap B|^{L} + |A|^{L_1}+|B|^{L_2}\\
&\leq |A|^L + |B|^L  + |A|^{L_1}+|B|^{L_2}.
\end{align*}
\end{proof}

\subsection{Reductive \texorpdfstring{$\boldsymbol\G$}{G}-Whitehead partitions}\label{sec:reductive}

Recall that a \GW partition $\Pbold$ of  a \rose\  $\sigma$ is   {\em reductive} if for some $v\in {\rm max}(P)$ the Whitehead automorphism $\varphi=(P,v)$  reduces the norm of $\sigma,$ i.e. $\|\sigma^\Pbold_v\|< \|\sigma\|$, and \emph{strongly reductive} if $(P,v)$ reduces the first coordinate of the norm, i.e. $\|\sigma^\Pbold_v\|_0< \|\sigma\|_0$.  

The strategy of our proof will require us to find reductive \GW partitions which are compatible with each other, so our next task is to determine where we can look for such partitions.  We first consider a pair of non-compatible partitions, and show how to find a partition which is compatible with both of them. 

In our definition of \GW partition, we did not allow $P$ to be a singleton.  For convenience, we now define a \emph{degenerate \GW partition} to be one of the form $\Pbold=(P,P^*,{\rm lk}(P))$ where $P=\{v\}$.  In this case,
the associated Whitehead automorphism $(P,v)$ is the inversion $\iota_v$ and $|\Pbold|_\sigma = |v|_\sigma$ for every $\sigma$.  In particular, a reductive \GW partition cannot be degenerate.

Suppose $\Pbold, \Qbold$ are \GW partitions which are not compatible, i.e. they do not commute and each of the sets 
$P \cap Q, P^*\cap Q, P\cap Q^*, P^*\cap Q^* $ is non-empty.  We will refer to these four intersections as \emph{quadrants}.
Two quadrants are \emph{opposite} if one is obtained from the other by switching sides of both  $\Pbold$ and $\Qbold$.

\begin{lemma}\label{opposite}  For any non-compatible partitions  $\Pbold, \Qbold$, there is a pair of opposite quadrants, each of which defines a (possibly degenerate) \GW partition with maximal element in $\{v^\pm,w^\pm\}$. 
\end{lemma}

\begin{proof} Let $v \in {\rm max}(P)$ and $w \in {\rm max}(Q)$.  Recall that $x \in {\rm double}(\Qbold)$ means that $x,x^{-1}$ both lie on the same side of $\Qbold$ and $x \in {\rm single}(\Qbold)$ means that $x, x^{-1}$ lie on opposite sides of $\Qbold$. 
We divide the proof into three cases.

Case 1:  $v \in {\rm double}(\Qbold)$ and $w \in {\rm double}(\Pbold)$.  
In this case, some quadrant contains an element of both $\{v^\pm\}$ and $\{w^\pm\}$.  Without loss of generality, we may assume that $v,w \in P \cap Q$.  We claim that, in this case, $(P \cap Q^*, w^{-1})$ is a \GW pair.

Let $C_v$ denote the component of $\G \setminus {\rm st}(w)$ which contains $v$.  Then $v \in Q$ implies $C_v \subset Q$.  Moreover, we have
\begin{enumerate}
\item if ${\rm lk}(x) \subseteq {\rm lk}(v)$ then either ${\rm lk}(x) \subseteq {\rm lk}(w)$ or $x \in C_v \subset Q$,
\item every component $C$ of $\G \setminus {\rm st}(w)$ with $C \neq C_v$ lies entirely in some component of $\G \setminus {\rm st}(v)$.
\end{enumerate}
The first property follows from the fact that if ${\rm lk}(x)\not\subseteq {\rm lk}(w)$ then $x$ is connected to $v$ via some vertex not in ${\rm lk}(w)$.  Hence $x$ and $v$ lie in the same component of $\G \setminus {\rm st}(w)$.  The second property follows from the fact that in order for ${\rm st}(v)$ to disconnect $C$, $C$ must intersect ${\rm st}(v)$ and hence it must contain $v$.  

We can now verify that $(P \cap Q^*, w^{-1})$ is \GW.  For if $x \in {\rm single}(P \cap Q^*)$ then either $x \in {\rm single}(Q^*)$, 
so ${\rm lk}(x) \subseteq {\rm lk}(w)$, or 
$x \in {\rm single}(P) \cap Q^*$, so ${\rm lk}(x) \subseteq {\rm lk}(v)$ and $x \notin Q$.  By (1), it follows that ${\rm lk}(x) \subseteq {\rm lk}(w)$.
If $x \in {\rm double}(P \cap Q^*) = double (P) \cap double (Q^*)$, then by (2), so is the component of $x$ in $\G \setminus {\rm st}(w)$.  
This proves that $(P \cap Q^*, w^{-1})$ is a \GW pair.  By symmetry, $(P^* \cap Q, v^{-1})$ is also a \GW pair.

Case 2:  $v \in {\rm double}(\Qbold), w \in {\rm single}(\Pbold)$.  In this case, $w,w^{-1}$ lie in opposite quadrants while $v, v^{-1}$ lie in adjacent quadrants.  It follows that some quadrant contains  an element of both $\{v,v^{-1}\}$ and $\{w,w^{-1}\}$.   Without loss of generality, we may assume that $v^{-1},w \in P^* \cap Q$.  We claim that $P^* \cap Q$ and $P \cap Q^*$ are \GW.

First consider $P \cap Q^*$.  Let $C_v$ be as above.  By assumption, $C_v \subset Q$.  Thus, the same argument as in case (1) applies to show that  $(P \cap Q^*, w^{-1})$ is a \GW pair. 

Next consider  $P^* \cap Q$.   Since $w \in {\rm single}(\Pbold)$, $w \leq v$ so any $x$ in ${\rm single}(P^* \cap Q)$  satisfies $x \leq v$.  
For ${\rm double}(P^* \cap Q)$, note that ${\rm lk}(w) \subseteq {\rm lk}(v)$ implies that every component of $\G \setminus {\rm st}(v)$ (other than the singleton $\{w\}$) is contained is some component of  $\G \setminus {\rm st}(w)$.  It follows that ${\rm double}(P^* \cap Q)$ is a union of components of $\G \setminus {\rm st}(v)$.   
Hence $(P^* \cap Q, v^{-1})$ is a \GW pair.

Case 3: $v \in {\rm single}(\Qbold), w \in {\rm single}(\Pbold)$.  This is only possible if ${\rm lk}(v)={\rm lk}(w)$. Since $v$ is a singleton in both partitions, $v$ and $v^{-1}$ lie in opposite quadrants.  Say $v \in P\cap Q$ and $v^{-1} \in P^* \cap Q^*$.  Then it is easy to see that $(P \cap Q, v)$ and $(P^* \cap Q^*, v^{-1})$ are \GW pairs.  Likewise, the opposite quadrants containing $w$ and $w^{-1}$ also give \GW pairs.
\end{proof}

We next need to add the condition that our \GW partitions be  reductive in certain situations. 
Let $\sigma$ be a \rose\  and  $\Pbold$ a \GW partition.  For the purpose of this discussion, we introduce a weaker notion of reductively:  we say that $\Pbold$ is \emph{0-reductive} for $\sigma$ if for some $v \in {\rm max}(P)$, 
$\|\sigma^\Pbold_v\|_0 \leq \|\sigma\|_0$.  

Fix $\sigma=(\SG,\alpha)$ and let $\mathcal W_0$ be a set of cyclically reduced words representing $\{ \alpha^{-1}(g) \mid g \in \mathcal G_0\}.$ Write $|\Pbold|_0 = \sum_{\word \in \mathcal W_0} |\Pbold|_\word$ and $|v|_0 = \sum_{\word \in \mathcal W_0} |v|_\word$.  Then $\Pbold$ is 
\begin{align*}
\textrm{0-reductive if\ } &\|\sigma^\Pbold_v\|_0 - \|\sigma\|_0= |\Pbold|_0 -|v|_0 \leq 0 \\
\textrm{reductive if\ } &\|\sigma^\Pbold_v\| - \|\sigma\| = |\Pbold|_\sigma -|v|_\sigma < \0  \in \Z \times \Z^{\mathcal G} \\
\textrm{strongly reductive if\ } &\|\sigma^\Pbold_v\|_0 - \|\sigma\|_0= |\Pbold|_0 -|w|_0 < 0.
 \end{align*}
In particular, strongly reductive $\Rightarrow$ reductive $\Rightarrow$ 0-reductive, but none of the converses hold.

\begin{lemma}\label{HLL} (Higgins-Lyndon Lemma)  Let $\sigma$ be a \rose\ and let $\Pbold$ and $\Qbold$ be non-compatible \GW partitions.
If  $\Pbold$ and $\Qbold$  are both $\sigma$-reductive 
then at least one of the quadrants $P\cap Q^*$, $P^*\cap Q$, $P\cap Q$ or $P^*\cap Q^*$ determines a 
$\sigma$-reductive \GW partition which is compatible with both $\Pbold$ and $\Qbold$.  If $\Pbold$ is strongly reductive and $\Qbold$ is 0-reductive, then one of the quadrants is strongly reductive.
\end{lemma}
\begin{proof}  
Let $\sigma=(\SG, \alpha)$.  By hypothesis, we can choose $v \in {\rm max}(P)$ and $w\in {\rm max}(Q)$ such that 
either  $|\Pbold|_{\sigma}-|v|_\sigma <\0$ and $ |\Qbold|_{\sigma}-|w|_\sigma <\0$ (case 1), 
or $|\Pbold|_0-|v|_0 < 0$ and $ |\Qbold|_0 - |w|_0 \leq  0 $ (case 2).

 Suppose first that there is exactly one quadrant which contains none of $\{v,v^{-1},w,w^{-1}\}$.  By changing sides of $\Pbold$ and $\Qbold$  if necessary, we may assume  this is $P \cap Q$.  Then both  $(P\cap Q^*,v)$ and $(P^*\cap Q,w)$ are \GW by Lemma~\ref{opposite}.  
  Since $P$ and $P^*\cap Q$ are disjoint and $v$ and $w$ don't commute, Lemma~\ref{commute} shows that $P\cap {\rm lk}(Q)=\emptyset$; similarly, $Q\cap {\rm lk}(P)=\emptyset$.  Thus $P$ and $Q$ are both in the complement of $L={\rm lk}(P)\cup {\rm lk}(Q)$. So by Lemma~\ref{ABL}, for every cyclically reduced word $\word$,
$$ |P\cap Q^*|^{{\rm lk}(v)}_\word + |P^*\cap Q|^{{\rm lk}(w)}_\word  \leq |P|^{{\rm lk}(v)}_\word+ |Q|^{{\rm lk}(w)}_\word $$

  Letting $\Pbold\cap \Qbold^*$ denote the \GW partition determined by $P\cap Q^*$ and $\Qbold\cap\Pbold^*$ the \GW partition determined by $Q\cap P^*,$ it follows  that
$$|\Pbold\cap \Qbold^*|_{\word}-|v|_\word + |\Pbold^*\cap \Qbold|_{\word}-|w|_\word\leq |\Pbold|_{\word}-|v|_\word + |\Qbold|_{\word}-|w|_\word .$$
In case 1, it follows that
$$|\Pbold\cap \Qbold^*|_{\sigma}-|v|_\sigma + |\Pbold^*\cap \Qbold|_{\sigma}-|w|_\sigma\leq |\Pbold|_{\sigma}-|v|_\sigma + |\Qbold|_{\sigma}-|w|_\sigma < \0,$$
so at least one of  $(P\cap Q^*,v)$ or $(P^*\cap Q,w)$ is reductive.  In case 2,
$$|\Pbold\cap \Qbold^*|_0-|v|_0+ |\Pbold^*\cap \Qbold|_0-|w|_0 \leq |\Pbold|_0-|v|_0 + |\Qbold|_0-|w|_0 < 0,$$
so one of $(P\cap Q^*,v)$ or $(P^*\cap Q,w)$ is strongly reductive.

Next suppose each quadrant contains an element  of $\{v, v^{-1},w,w^{-1}\},$  and say $v\in P\cap Q$.  This forces $w\in Q\cap P^*, w^{-1}\in P\cap Q^*$ and $v^{-1}\in P^*\cap Q^*$.  Since $v$ is a singleton in $Q$, and $w$ is a singleton in $P$, we have ${\rm lk}(v)={\rm lk}(w)=L$ and by Lemma~\ref{opposite}, all four quadrants are \GW.  Recall that 
$|P|^L_\word =|P^*|^L_\word$, so applying Lemma~\ref{ABL} to both pairs of opposite quadrants gives
\begin{align*}
&(|P\cap Q^*|^L_\word-|w|_\word+ |P^*\cap Q|^L_\word-|w|_\word)+(|P^*\cap Q^*|^L_\word-|v|_\word+ |P\cap Q|^L_\word-|v|_\word) \\
&\leq (|P|^L_\word+ |Q|^L_\word-2|w|_\word) + ( |P^*|^L_\word + |Q|^L_\word-2|v|_\word)\\
& = 2(|P|^L_\word -|v|_\word) + 2(|Q|^L_\word - |w|_\word)
\end{align*}
In case 1 we obtain
\begin{center}
$(|\Pbold\cap \Qbold^*|_\sigma-|w|_\sigma)+ (|\Pbold^*\cap \Qbold|_\sigma-|w|_\sigma)+(|\Pbold^*\cap \Qbold^*|_\sigma-|v|_\sigma)+ (|\Pbold\cap \Qbold|_\sigma-|v|_\sigma)$ \\
$\leq 2(|\Pbold|_\sigma -|v|_\sigma) + 2(|\Qbold|_\sigma - |w|_\sigma) < \0$
\end{center}
so at least one of the quadrants is  reductive.  In case 2 we have
\begin{center}
$(|\Pbold\cap \Qbold^*|_0-|w|_0)+ (|\Pbold^*\cap \Qbold|_0-|w|_0)+(|\Pbold^*\cap \Qbold^*|_0-|v|_0)+ (|\Pbold\cap \Qbold|_0-|v|_0)$ \\
$\leq 2(|\Pbold|_0 -|v|_0) + 2(|\Qbold|_0 - |w|_0) < 0$
\end{center}
so one of the quadrants is strongly reductive.

The remaining possibility is that only two quadrants contain elements of $\{v, v^{-1},w,w^{-1}\}$.  In this case, we may assume $v, w\in P\cap Q$ and $v^{-1},w^{-1}\in P^*\cap Q^*$.  Here again ${\rm lk}(v)={\rm lk}(w)=L$, and   $(P\cap Q,v)$ and $(P^*\cap Q^*,w^{-1})$ are both \GW by Lemma~\ref{opposite}.  Applying Lemma~\ref{ABL} gives
$$ (|\Pbold\cap \Qbold|_\word-|v|_\word) + (|\Pbold^*\cap \Qbold^*|_\word-|w|_\word) \leq (|\Pbold|_\word-|v|_\word) + 
(|\Qbold|_\word-|w|_\word)$$
and arguing as above we conclude that one of these quadrants is reductive (case 1)  or strongly reductive (case 2).  

Finally, note that the requirement that the chosen quadrant define a partition compatible with both $\Pbold$ and $\Qbold$ is immediate from the fact that every quadrant is contained in one side of $\Pbold$ and one side of $\Qbold$.
\end{proof}

We have shown that any two marked Salvettis, $\sigma,  \sigma',$  can be joined by a path in $\KG$ consisting of a sequence of Whitehead moves.  We call such a path a \GW path.

\begin{theorem}\label{peak}(Peak Reduction).  Let $(P,v)$ and $(Q,w)$ be two reductive \GW moves from  $\sigma$. Then there is a \GW path from $\sigma^\Pbold_v$ to $\sigma^\Qbold_w$  which passes only through \roses\ $\tau$   with $\|\tau\|<\|\sigma\|$. 
\end{theorem}

\begin{proof}  First observe that in the case where $\Pbold, \Qbold$ are compatible and $v = w$,  it follows from Remark \ref{maxPQ} that  $\sigma^\Pbold_v$ and $\sigma^\Qbold_w$ differ by a single \GW move so there is nothing to prove.  

Consider the case where $\Pbold, \Qbold$ are compatible and $v \neq w$.  If the edges $e_v, e_w$ in $S^{\Pbold, \Qbold}$ do not join the same two vertices, then the hyperplanes $\Hcal= \{H_v,H_w\}$ through these edges form a tree-like set in $S^{\Pbold, \Qbold}$. In this case, setting $\tau=\sigma^{\Pbold,\Qbold}_\Hcal$, we obtain a \GW path
$\sigma^\Pbold_v \to  \tau \to   \sigma^{\Qbold}_w$.  Since $(P,v)$ and $(Q,W)$ are both reductive, Lemma \ref{norm2} gives
$$\|\tau\| =\|\sigma\| + (|\Pbold|_\sigma - |v|_\sigma) + (|\Qbold|_\sigma - |w|_\sigma) < \|\sigma\|.$$

The only situation in which $e_v,e_w$ can join the same pair of vertices is if  $v,w$ are singletons in both partitions, say $v,w \in P \cap Q$ and $v^{-1},w^{-1} \in P^* \cap Q^*$.  In this case, 
$(P,w)$ and $(Q, v)$ are also \GW pairs.  Suppose $|v|_\sigma \leq |w|_\sigma$.  Then by Lemma \ref{norm2},
$$\|\sigma^\Pbold_w\| = \|\sigma\| + (|\Pbold|_\sigma - |w|_\sigma) \leq \|\sigma\| + (|\Pbold|_\sigma - |v|_\sigma)
 < \|\sigma\| $$
so $\sigma^\Pbold_v \to  \sigma^\Pbold_w  \to   \sigma^{\Qbold}_w$ is the desired path.  
If $|v|_\sigma > |w|_\sigma$,  use $\sigma^\Pbold_v \to  \sigma^\Qbold_v  \to   \sigma^{\Qbold}_w$ instead.

Next, suppose $\Pbold$ and $\Qbold$ are not compatible.
Apply the Higgins-Lyndon lemma to find $\Rbold$ compatible with both $\Pbold$ and $\Qbold$, with $\Rbold$ reductive,
i.e., $\|\sigma^\Rbold_u\|  < \|\sigma\|$.   Then by the discussion above, there are \GW paths  from $\sigma^\Pbold_v$ to $\sigma^{\Rbold}_{u}$ and from $\sigma^{\Rbold}_{u}$ to $\sigma^{\Qbold}_w$ satisfying the required condition. 
\end{proof}

\begin{theorem}\label{s-peak}(Strong Peak Reduction).  Let $(P,v)$ and $(Q,w)$ be two \GW partitions such that 
$\|\sigma^\Pbold_v\|_0 < \|\sigma\|_0$ and $\|\sigma^\Qbold_v\|_0 \leq \|\sigma\|_0$. Then there is a \GW path from $\sigma^\Pbold_v$ to $\sigma^\Qbold_w$  which passes only through \roses\ $\tau$   with $\|\tau\|_0<\|\sigma\|_0$. 
\end{theorem}

\begin{proof}  Let $\mathcal W_0$ be a set of cyclically reduced words representing $\{ \alpha^{-1}(g) \mid g \in \mathcal G_0\}. $  Write $|\Pbold|_0 = \sum_{\word \in \mathcal W_0} |\Pbold|_\word$ and $|v|_0 = \sum_{\word \in \mathcal W_0} |v|_\word$.  Then 
\begin{align*}
\|\sigma^\Pbold_v\|_0 - \|\sigma\|_0 = |\Pbold|_0 -|v|_0 &<0 \\
 \|\sigma^\Qbold_v\|_0 - \|\sigma\|_0= |\Qbold|_0 -|w|_0 &\leq 0.
 \end{align*}
 We now proceed as in the proof of the previous theorem.  In the case where $\Pbold, \Qbold$ are compatible and $e_v, e_w$ join different vertices in $S^{\Pbold, \Qbold}$, set  $\tau=\sigma^{\Pbold,\Qbold}_\Hcal$ and note that 
 $$\|\tau\|_0 =\|\sigma\|_0 + (|\Pbold|_0 - |v|_0) + (|\Qbold|_0 - |w|_0) < \|\sigma\|_0.$$
If $e_v, e_w$  join the same vertices and $|v|_0 \leq |w|_0$, then
$$\|\sigma^\Pbold_w\|_0 = \|\sigma\|_0 + (|\Pbold|_0 - |w|_\0) \leq \|\sigma\|_0 + (|\Pbold|_0 - |v|_0) < \|\sigma\|_0$$
If $e_v, e_w$  join the same vertices and $|v|_0 > |w|_0$, then
$$\|\sigma^\Qbold_v\|_0 = \|\sigma\|_0 + (|\Qbold|_0 - |v|_\0) < \|\sigma\|_0 + (|\Qbold|_0 - |w|_0) \leq  \|\sigma\|_0.$$
Hence either  $\sigma^\Pbold_v \to  \sigma^\Pbold_w  \to   \sigma^{\Qbold}_w,$ or  $\sigma^\Pbold_v \to  \sigma^\Qbold_v  \to   \sigma^{\Qbold}_w,$ gives the desired path.

\end{proof}

\begin{corollary}\label{existence} If $\|\sigma\|$ is not minimal, then there is a  strongly reductive \GW move from $\sigma$.    
\end{corollary}
\begin{proof}  As observed in the proof of Proposition \ref{connected}, there is a path, 
$\sigma = \sigma_1 , \sigma_2 , \dots \sigma_k=(\SG, id),$
of \GW moves from $\sigma$ to the unique \rose\  with minimal 0-norm, $\|(\SG,id) \|_0$.    Consider the sequence of 0-norms,  $\|\sigma\|_0 , \|\sigma_2\|_0 , \dots \|\sigma_k\|_0$.   Using Theorem \ref{s-peak}, we can reduce peaks in this sequence to obtain a \GW path $\sigma=\tau_0, \tau_1,\dots ,\tau_j=(\SG,id)$ which begins downward, that is, with $\|\sigma\|_0 >\|\tau_1\|_0$ 
\end{proof}

\begin{corollary}\label{finite}  Let $N_0=\|(\SG,id)\|_0$.  For any $N \geq N_0$, there are  finitely many \roses\ $\sigma$ with $\|\sigma\|_0 \leq N$.
\end{corollary}
\begin{proof}  We observed in Lemma \ref{unique} that  $(\SG,id)$ is the unique \rose\ with minimal 0-norm.  By Corollary \ref{existence}, if  $\|\sigma\|_0 \leq N$, there is a \GW path of length at most $N-N_0$ to $(\SG, id)$.  Since the number of Whitehead moves at any \rose\ is bounded, the Corollary follows.
\end{proof}

\begin{proposition} The set of  \roses\  is well-ordered with respect to the norm $\|\ \ \|$.
\end{proposition}
\begin{proof}
This follows immediately from Corollary \ref{finite}.  Let $N=\| \sigma \|_0$.  Since $\|\tau\| < \|\sigma\|$ implies $\|\tau\|_0 \leq \|\sigma\|_0$, there are only finitely many such $\tau$.  Hence there can be no infinite decreasing chain of \roses.  
\end{proof}

For the proof of contractibility, we will also need the following stronger form of Lemma~\ref{HLL}.

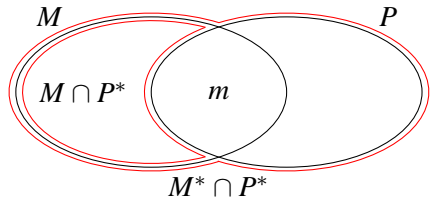
\begin{figure}
\begin{center}
  \begin{tikzpicture} [xscale=.45, yscale=.25]
    \draw (0,0) circle (4cm); 
  \node (M) at (-3,4) {$M$};
   \draw [red] (65:3.8) arc (65:295:3.8cm);
   \draw [red] (61.5:4.2) arc (61.5:298.5:4.2cm);
  \begin{scope}[xshift=4cm]
   \draw (0,0) circle (4cm);
   \node (P) at (3,4) {$P$}; 
     \draw [red] (118.5:4.2) arc (118.5:-119:4.2cm);
       \draw [red] (125:4.2) arc (125:235:4.2cm);
    \end{scope}
    
    \node (m) at (2,0) {$m$};
     \node (MPstar) at (-2,0) {$M\cap P^*$};
      \node (MP) at (2,-5) {$M^*\cap P^*$};
   \end{tikzpicture}
   \caption{Setup for the Pushing Lemma}\label{pushinglemma}
\end{center}
\end{figure}

\begin{lemma}\label{pushing}  (Pushing Lemma)
Fix a \rose\  $\sigma$. Suppose that  $(M,m)$ is a reductive \GW pair such that at $\sigma$,
\begin{enumerate}
\item ${\rm lk}(\Mhat)$ is maximal among links of reductive \GW partitions, and 
\item  $(M, m)$  is maximally reductive among \GW pairs $(Q,w)$ with ${\rm lk}(\Qbold)={\rm lk}(\Mhat)$. 
\end{enumerate}
 Let  $\Pbold$ be a reductive \GW partition that is not compatible with $\Mhat$.  Choose sides so that $m\in M \cap P$ .  Then at least one of $M\cap P^*$ or $M^*\cap P^*$ determines  a reductive \GW partition whose link is equal to  ${\rm lk}(\Pbold)$  (see Figure~\ref{pushinglemma}).
\end{lemma}

\begin{proof}  
$\Pbold$ is reductive, so for some $v \in {\rm max}(\Pbold)$, $|P|_\sigma - |v|_\sigma < \0$.  
We will apply Lemma~\ref{opposite} to $\Pbold,\Mhat$.  

We first consider the case where  $m \in {\rm single}(\Pbold)$ and $v \in {\rm single}(\Mhat)$.  This corresponds to case (3) in the proof of Lemma~\ref{opposite}.   In this case, ${\rm lk}(v)={\rm lk}(m)$ and the lemma produces opposite quadrants $X,Y$ such that $(X,v)$ and $(Y,v^{-1})$ are \GW pairs, as well as and opposite quadrants $X',Y'$ such that $(X',m)$ and $(Y',m^{-1})$ are \GW pairs.  
Applying Lemma~\ref{ABL} to both pairs of opposite quadrants gives
$$|{\Xbold}|_\sigma + |\Ybold|_\sigma +|{\Xbold'}|_\sigma + |\Ybold'|_\sigma \leq 2(|\Pbold|_\sigma + |\Mhat|_\sigma)$$
Hence 
\begin{align*}
&(|{\Xbold}|_\sigma -|v|_\sigma) + (|\Ybold|_\sigma -|v|_\sigma) +( |{\Xbold'}|_\sigma -|m|_\sigma) +( |\Ybold'|_\sigma -|m|_\sigma) \\
& \leq  2(|\Pbold|_\sigma -|v|_\sigma)  + 2( |\Mhat|_\sigma - |m|_\sigma)\\
& < 2( |\Mhat|_\sigma - |m|_\sigma)
\end{align*}
Two of these quadrants, say $Y$ and $Y',$ lie in $P^*$.   By hypothesis (2), $(M,m)$ is maximally reductive, so 
$$2(|\Mhat|_\sigma -|m|_\sigma) \leq( |{\Xbold}|_\sigma -|v|_\sigma) + ( |{\Xbold'}|_\sigma -|m|_\sigma)$$
and we conclude that  
$$(|{\Ybold}|_\sigma -|v|_\sigma) + (|\Ybold'|_\sigma -|m|_\sigma) < \0. $$
Thus one of the pairs $(Y,v^{-1})$ or $(Y',m^{-1})$ is reductive and satisfies the requirements of the lemma.

In all other cases,  Lemma \ref{opposite} gives
 an opposite pair of quadrants $X$ and $Y$,  with maximal elements $x\in \{v^\pm\}, y \in  \{m^{\pm}\}$, that define (possibly degenerate) \GW partitions.   By Lemma \ref{ABL} we have
$$(|{\Xbold}|_\sigma -|v|_\sigma) + (|\Ybold|_\sigma -|m|_\sigma) \leq (|\Pbold|_\sigma -|v|_\sigma) + (|\Mhat|_\sigma -|m|_\sigma) < (|\Mhat|_\sigma -|m|_\sigma).$$
By hypothesis (2), $ |\Mhat|_\sigma -|m|_\sigma \leq |{\Ybold}|_\sigma -|m|_\sigma$, so we conclude that
 $(X,x)$ is reductive.
 If  $X \subset P^*$, we are done.   If $X \subset P$, then $Y \subset P^*,$ so $P^*$ contains $y \in \{m^\pm\}$ and we must have $m \in {\rm single}(P)$.  This implies that ${\rm lk}(m) \subseteq {\rm lk}(v)$, so by hypothesis (1),  ${\rm lk}(m)={\rm lk}(v)$.
It then follows from hypothesis (2) that $ |\Mhat|_\sigma -|m|_\sigma \leq |{\Xbold}|_\sigma -|v|_\sigma$, and we conclude that $(Y,y)$ is also reductive.
\end{proof}

\subsection{Contractibility of the  \texorpdfstring{$\boldsymbol\G$}{G}-spine \texorpdfstring{$\KG$}{KG}}

In this section we prove our main theorem.

\begin{theorem}  \label{contract} For any right-angled Artin group $\AG$, the  $\G$-spine $\KG$ is contractible.
\end{theorem}

The proof will make frequent use the following lemma, which is standard in the topology of posets, and dates back to work of Quillen \cite{Quillen}. 
 
 \begin{lemma}\label{poset} (Poset Lemma)  Let $X$ be a poset and $f\colon X\to X$   a poset map with the property that $ x\leq f(x)$ for all $x\in X$ (or $x\geq f(x)$ for all $x\in X$). Then f induces a deformation retraction from the geometric realization  of $X$ to the geometric realization of the image $f(X).$
 \end{lemma}

\begin{proof}[Proof of Theorem \ref{contract}]    We view the spine $\KG$ as the union of stars of \roses.  By Lemma~\ref{unique} there is a unique  \rose\  $(\SG,id)$ of minimal norm, and we start with its (contractible) star.  We build the entire spine by gluing on stars of \roses\  in increasing order.

When we add a \rose, we need to check that we are gluing along something contractible.  So fix a \rose\  $\sigma$, and let $K_{<\sigma}$ be the union of stars of \roses\  $\tau$ with $\|\tau\| < \|\sigma\|$.  The intersection  ${\rm st}(\sigma)\cap K_{<\sigma}$ consists of marked blowups $\sigma^\PP$ which can be collapsed to a \rose\  of smaller norm;
 here  $\PP=\{\Pbold_1, \dots \Pbold_k\}$ is a set of compatible \GW partitions, which we will refer to as an ``ideal forest."   We can identify ${\rm st}(\sigma)\cap K_{<\sigma}$ with the geometric realization of the poset of such ideal forests, ordered by inclusion.   To prove that ${\rm st}(\sigma)\cap K_{<\sigma}$ is contractible, we will repeatedly apply the Poset Lemma to retract this poset to a single point.

First note that by Corollary~\ref{existence},  ${\rm st}(\sigma)\cap K_{<\sigma}$ is non-empty.  
 By Corollary~\ref{reductive}, if $\PP$ is in ${\rm st}(\sigma)\cap K_{<\sigma}$, then $\PP$ contains a $\sigma$-reductive \GW partition.  Therefore, the map that throws out the non-reductive \GW partitions from each $\PP$ in ${\rm st}(\sigma)\cap K_{<\sigma}$ is a poset map, and we can use  the Poset Lemma to retract ${\rm st}(\sigma)\cap K_{<\sigma}$ to  its image, which is the subposet  $\mathcal R$ spanned by ideal forests all of whose \GW partitions are reductive.  
 
 Now choose a reductive pair $(M,m)$ satisfying the maximality conditions of Lemma \ref{pushing}.   We  will ultimately retract $\mathcal R$ to the ideal forest consisting of the single partition  $\{\Mhat\}$.   

If all of the ideal forests  in $\mathcal R$ are compatible with $\Mhat$, then we can contract $\mathcal R$ to $\{\Mhat\}$ via the poset maps $\PP \to\PP\cup\{\Mhat\} \to \{\Mhat\}$.
If not, choose a reductive \GW partition $\Pbold$  such that
\begin{enumerate}
\item $\Pbold$ and $\Mhat$ are not compatible, and
\item  the side $P$ containing $m$ is maximal among sides of all such partitions, i.e., if  $\Qbold$ is a reductive \GW partition with $P \subset Q$, then $\Qbold$ is compatible with $\Mhat$.  
\end{enumerate}
Note that if $m \in P \subset Q$, then $m \notin {\rm lk}(\Qbold)$, so $\Mhat$ and $\Qbold$ do not commute. It follows that the only way they can be compatible is if $M \subset Q$.  Thus, the second condition can be restated as,
 $(2')$ if $P \subset Q$, then $M \subset Q$.

By the Pushing Lemma, one of $M\cap P^*$ or $M^*\cap P^*$ determines a reductive \GW partition whose link is equal to ${\rm lk}(\Pbold)$, call it $\Pbold_0$. 
We claim that 
$$\PP\mapsto \begin{cases} \PP \cup\{\Pbold_0\}   & \mbox{if } \Pbold\in\PP \\  \PP  & \mbox{if }  \Pbold\not\in\PP \end{cases}$$   is a poset map from $\mathcal R$ to itself.    If $\Pbold\in\PP$ then $\Qbold\in \PP$ implies that $\Qbold$ is compatible with $\Pbold$, so we have to check that any such  $\Qbold$ is also compatible with $\Pbold_0$.  

If $\Qbold$ commutes with $\Pbold$, then it also commutes with $\Pbold_0$ since they have the same link.  Otherwise,
compatibility implies that  some side $Q$ of $\Qbold$ is either contained in $P$ or contains $P$.  If  $Q \subset P$, then  $Q \cap P^*=\emptyset$, so  $\Qbold$ is compatible with both $M\cap P^*$ and $M^*\cap P^*$.    If $P \subset Q$, then  by  $(2')$, $M \cup P  \subset Q$.  It follows that   $M\cap P^* \subset Q$ and $Q^* \subset (M^* \cap P^*)$,  so both of these quadrants are again compatible with $\Qbold$.  This proves the claim.

This map clearly satisfies the hypotheses of the Poset lemma, so $\mathcal R$ retracts to the image, in which every ideal forest that contains $\Pbold$ also contains $\Pbold_0$. Now map this image to itself by the poset map which throws $\Pbold$ out of every $\PP$ that contains it. The Poset Lemma applies again, and the image is now the subcomplex of $\mathcal R$ spanned by all reductive ideal forests which do not contain $\Pbold$.

 Repeat this process until  every \GW partition that is not compatible with $\Mhat$ has been eliminated.  The resulting poset can be retracted to the single point $\{\Mhat\}$ as described above.     
\end{proof} 

\section{Outer space for  \texorpdfstring{$\Out$}{Out(Fn)}}

In this final section we briefly discuss the problem of finding an outer space for the entire group $\Out$.  If there are no vertices $v,w$ in $\G$ with ${\rm st}(v)\subseteq {\rm st}(w)$, then $\UAG=\Out$ and the spine $\KG$ is quasi-isometric to $\Out$, so nothing more needs to be done. 

 In particular, if $\AG$ is a free group then $\KG$ is the spine of Culler-Vogtmann's outer space $\mathcal O_n$.  The complete space $\mathcal O_n$ is obtained by  giving edges of marked graphs positive real lengths, then taking projective classes of the resulting metric graphs.  The spine $K_n$ naturally embeds into $\mathcal O_n,$ by sending  vertices of $K_n$  to metric graphs with all edges of equal length.  
We can do the same thing for general $\KG$ by considering the cubes in a $\G$-complex to be arbitrary rectilinear parallelepipeds, i.e. giving each set of edges dual to a hyperplane the same positive real length, then taking the projective class of the resulting metric space.  We call the space we obtain in this way $\Sigma_\G$, and $\KG$ embeds into it as in the free group case, as an equivariant deformation retract.  

The opposite extreme from the free group is the free abelian group $\Z^n$.  In this case  $\Out\iso {\rm GL}(n,\Z)$  acts properly on the symmetric space ${\rm SL}(n,\R)/{\rm SO}(n)$, which we can regard   as a space of projective classes of flat tori.  For general $\G$,  abelianization $\AG\to\Z^n$ induces a map $\Out\to {\rm GL}(n,\Z)$ which sends the twist subgroup $T_\G$ injectively into ${\rm SL}(n,\Z)$.  There is a natural contractible subspace $\mathbb D_\G$  of the symmetric space on which $T_\G$ acts; this space $\mathbb D_\G$ can be regarded as a space of projective classes of certain ``flat Salvettis;" where we regard the Salvetti as a subcomplex of a flat $n$-torus and take the induced path metric.  This gives us a restricted class of allowable metrics on Salvetti complexes, which can be extended to a restricted class of allowable metrics on  $\G$-complexes.  Cubes are no longer necessarily rectilinear; some are allowed to become parallelepipeds.

We can then define an outer space $\OG$ on which all of $\Out$ acts as follows:
A point in $\OG$ is an equivalence class of triples $(X,d,h)$ where
\begin{itemize}
\item $X$ is a $\G$-complex
\item $d$ is an allowable metric on $X$
\item $h\colon X\to \SG$ is an arbitrary homotopy equivalence
\item $(X,d,h) \sim (X',d',h')$ if there is an isometry $f\colon X\to X'$ with $h'\circ f\simeq h$.
\end{itemize}
The action of $\varphi\in \Out$ on $\OG$ is by composition:  realize $\varphi$ by $f\colon\SG\to\SG$, then $\varphi(X,d,h)=(X,d,f\circ h)$.

Both $\Sigma_\G$ and $\mathbb D_\G$ embed into this space $\OG$. In many cases it is easy to reduce the contractibility of $\OG$ to the contractility of $\Sigma_\G$, proved in this paper.  However the general case seems to be quite subtle and not at all short, so this will form the contents of a second paper.

%
%
\bibliographystyle{gtart}

\end{document}